\documentstyle[12pt]{article}

\newtheorem{LE}{Lemma}
\newtheorem{TH}{Theorem}
\newtheorem{PR}{Proposition}
\newtheorem{CO}{Corollary}
\newtheorem{EX}{Example}
\newcommand{\lh}{\leftharpoonup}
\newcommand{\rh}{\rightharpoonup}
\newcommand{\pf}{\medskip\noindent{\sc Proof: }}
\newcommand{\la}{\longrightarrow}
\newcommand{\rTo}{\,\longrightarrow\,}
\newcommand{\qed}{$\Box$}
\newcommand{\SL}{\sl}

\title{Oriented Quantum Algebras and Invariants of Knots and Links} 
\author{Louis H. Kauffman\thanks{Research supported in part by NSF Grant
DMS 920-5227} $\;\;$
and $\;\;$
David E. Radford\thanks{Research supported in part by NSF Grant
DMS 980 2178}\\  \\ 
Department of Mathematics, Statistics \\
and Computer Science (m/c 249)    \\
851 South Morgan Street   \\
University of Illinois at Chicago\\
Chicago, Illinois 60607-7045}
\begin{document}
\maketitle
\begin{abstract}
{\small \rm
In 
\cite{KRO}
oriented quantum algebras were motivated and introduced  in a very natural categorical setting within the context of  knots and links, some examples were discussed and a rudimentary theory of oriented quantum algebras was sketched. Invariants of knots and links can be computed from oriented quantum algebras.

Here we continue the study of oriented quantum algebras. We view them from a more algebraic perspective and develop a more detailed theory for them and their associated invariants. We study a class of examples associated with the HOMFLY polynomial in depth.
}
\end{abstract}
\setcounter{section}{-1}
\section{Introduction}
In 
\cite{KRO}
the notion of oriented quantum algebra was introduced and motivated from a topological point of view. Oriented quantum algebras give rise to regular isotopic invariants of oriented $1$--$1$ tangles and twist oriented quantum algebras, which are oriented quantum algebras with additional structure, give rise to regular isotopy invariants of oriented knots and links. Twist oriented quantum algebras are important in that their associated invariants include nearly all of the quantum links invariants know at present. 

In this paper we develop a general theory of oriented quantum algebras, twist oriented quantum algebras and their resulting invariants. Quantum algebras have oriented quantum algebra structures and twist quantum algebras have twist oriented quantum algebra structures. Quasitringular Hopf algebras are examples of quantum algebras and ribbon Hopf algebras are examples of twist oriented quantum algebras. Not every oriented quantum algebra is accounted for by a quantum algebra. There are examples of twist oriented quantum algebras on the algebra ${\rm M}_n(k)$ of $n {\times} n$ matrices over a field $k$ which we explore in depth related to the HOMFLY polynomial. 

This paper is the second in a series of four on the theory of oriented quantum algebras and related structures. In the first paper \cite{KRO} we motivate the definition of oriented quantum
algebras by considering the structure of a functor from the tangle category to a category associated with an
algebra. The first paper shows how all known quantum link invariants arise from oriented quantum algebras.
In the third 
\cite{KRNexT}
we introduce the notion of oriented quantum coalgebra and related concepts. The notion of oriented quantum coalgebras are a bit more general than the dual of oriented quantum algebra. A theory of oriented quantum algebras is developed along the lines of the theory of quantum coalgebras found in 
\cite{QC}.
In the fourth
\cite{KRNexTNexT}
paper of this series the connection between the state sum description of many quantum link invariants and their formulation in terms of oriented quantum algebras (and coalgebras) is discussed in great detail. 

Throughout $k$ is a field and all vector spaces are over $k$. The set of non-zero elements of $k$ is denoted by $k^\star$ and all vector spaces over $k$.
\section{Preliminaries}\label{secPre}
For vector spaces $U$ and $V$ over $k$ we will denote the tensor product $U{\otimes}_kV$ by $U {\otimes} V$ and the identity map of $V$ by $1_V$. If $t$ is a linear endomorphism of $V$ then an element $v \in V$ is $T$-{\em invariant} if $t(v) = v$. If $A$ is an algebra over $k$ we shall let $1_A$ also denote the unit of $k$. Then meaning $1_V$ should always be clear from context.

Let $A$ be an algebra over the field $k$ and let $\rho = \sum_{\imath = 1}^r a_\imath {\otimes} b_\imath \in A{\otimes} A$. We set 
$$
\rho_{1\, 2} = \sum_{\imath = 1}^r a_\imath {\otimes} b_\imath {\otimes} 1, \;\; 
\rho_{1\, 3} = \sum_{\imath = 1}^r a_\imath {\otimes} 1 {\otimes} b_\imath \;\; 
\mbox{and} \;\;
\rho_{2\, 3} = \sum_{\imath = 1}^r 1 {\otimes} a_\imath {\otimes} b_\imath.
$$
The quantum Yang--Baxter equation for $\rho$ is $\rho_{1\, 2} \rho_{1\, 3} \rho_{2\, 3}  = \rho_{2\, 3} \rho_{1\, 3}\rho_{1\, 2}$. There is an important class of solutions to the quantum Yang--Baxter equation for the algebra $A = {\rm M}_n(k)$ of all $n{\times}n$ matrices over $k$. For $1 \leq \imath, \jmath \leq n$ let $E_{\imath \, \jmath} \in {\rm M}_n(k)$ be the $n{\times}n$ matrix which has a single non-zero entry, the value $1$ located in the $i^{th}$ row and $\jmath^{th}$ column. Then $\{ E_{\imath \,\jmath} \}_{1 \leq \imath, \jmath \leq n}$ is the standard basis for ${\rm M}_n(k)$ and $E_{\imath \, \jmath} E_{\ell \, m} = \delta_{\jmath \, \ell} E_{\imath \, m}$ for all $1 \leq \imath, \jmath, \ell, m \leq n$.
\begin{EX}\label{ExamRaBC}
Let $n \geq 2$, $a, b\!c \in k^\star$ satisfy $a^2 \neq b\!c, 1$ and let 
$$
{\sf B} = \{ b_{\imath\,  \jmath}\,|\, 1 \leq \imath < \jmath \leq n\},\qquad {\sf C} = \{ c_{\jmath \, \imath}\, 1 \leq \imath < \jmath \leq n\}
$$
be indexed subsets of $k^\star$ such that $b_{\imath \, \jmath}c_{\jmath \, \imath} = b\!c$ for all $1 \leq \imath < \jmath \leq n$. Then 
$$
\rho_{a, {\sf B}, {\sf C}} = \sum_{1 \leq \imath < \jmath \leq n}(a - \frac{b\!c}{a})E_{\imath \, \jmath} {\otimes} E_{\jmath \, \imath} + \sum_{\imath = 1}^naE_{\imath \, \imath} {\otimes} E_{\imath \, \imath} + 
 \sum_{1 \leq \imath < \jmath \leq n}(b_{\imath\, \jmath}E_{\imath \, \imath} {\otimes} E_{\jmath \, \jmath}  + c_{\jmath \, \imath}E_{\jmath \, \jmath} {\otimes} E_{\imath \, \imath}) 
$$
satisfies the quantum Yang--Baxter equation.
\end{EX}
That  $\rho_{a, {\sf B}, {\sf C}}$ satisfies the quantum Yang--Baxter equation follows by
\cite[Lemma 4 and (37)]{RBracket}.

By $A^{op}$ we mean the $k$-algebra whose underlying vector space is $A$ and whose multiplication is given by $a{\cdot}b = ba$ for all $a, b \in A$. An {\em oriented quantum algebra} over the field $k$ is a quadruple $(A, \rho, t_{\sf d}, t_{\sf u})$, which we sometimes informally designate by $A$, where $A$ is an algebra over $k$, $\rho \in A {\otimes}A$ is invertible and $t_{\sf d}, t_{\sf u}$ are commuting algebra automorphisms of $A$, such that 
\begin{flushleft}
{\rm (qa.1)} 
$(1_A{\otimes} t_{\sf u})(\rho)$ and $(t_{\sf d} {\otimes}1_A)(\rho^{-1})$ are inverses in $A{\otimes}A^{op}$,
\vskip1\jot
{\rm (qa.2)} 
$\rho = (t_{\sf d} \otimes t_{\sf d})(\rho) = (t_{\sf u} \otimes t_{\sf u})(\rho)$ and
\vskip1\jot
{\rm (qa.3)}
$\rho_{1\, 2} \rho_{1\, 3} \rho_{2\, 3}  = \rho_{2\, 3} \rho_{1\, 3}\rho_{1\, 2}.$
\end{flushleft}
In 
\cite{KRO}
the algebra automorphisms $t_{\sf d}$ and $t_{\sf u}$ are denoted by $D$ and $U$ which are more suitable labels for a diagrammatic treatment of oriented quantum algebras.

An oriented quantum algebra $(A, \rho, t_{\sf d}, t_{\sf u})$ over $k$ is {\em standard} if $t_{\sf d} = 1_A$ and is {\em balanced} if $t_{\sf d}  = t_{\sf u}$,  in which case we write $(A, \rho, t)$ for $(A, \rho, t_{\sf d}, t_{\sf u})$, where $t = t_{\sf d} = t_{\sf u}$. Any oriented quantum algebra $(A, \rho, t_{\sf d}, t_{\sf u})$ over $k$ gives rise to a standard one.  For by applying the algebra automorphisms $t_{\sf u} {\otimes} 1_A$ and $1_A {\otimes} t_{\sf d}$ of $A {\otimes} A^{op}$ to both sides of the equations of (qa.1)  one sees that:
\begin{PR}\label{ToStandard}
If $(A, \rho, t_{\sf d}, t_{\sf u})$ is an oriented quantum algebra over $k$ then $(A, \rho, t_{\sf u} {\circ} t_{\sf d}, 1_A)$ and $(A, \rho, 1_A, t_{\sf d}{\circ} t_{\sf u})$ are also. \qed
\end{PR}
\medskip

\noindent
The oriented quantum algebra $(A, \rho, 1_A, t_{\sf d}{\circ} t_{\sf u})$ is the {\em standard oriented quantum algebra associated with} $(A, \rho, t_{\sf d}, t_{\sf u})$. 

Important examples of  balanced oriented quantum algebras arise from Example \ref{ExamRaBC}.
\begin{EX}\label{ExamMnkRaBC}
Let $n \geq 2$, $a, b\!c \in k^\star$ satisfy $a^2 \neq b\!c, 1$ and suppose $\omega_1, \ldots, \omega_n \in k^\star$ satisfy 
$$
\omega_\imath^2 = \left(\frac{a^2}{b\!c}\right)^{\imath-1}\omega_1^2
$$
for all $1 \leq \imath \leq n$. Then $({\rm M}_n(k), \rho_{a, {\sf B}, {\sf C}}, t)$ is a balanced oriented quantum algebra, where 
$$
t(E_{\imath \, \jmath}) = \left(\frac{\omega_\imath}{\omega_\jmath}\right)E_{\imath \, \jmath}
$$
for all $1 \leq \imath, \jmath \leq n$ and $\rho_{a, {\sf B}, {\sf C}}$ is described in Example \ref{ExamRaBC}.
\end{EX}

The assertions of Example \ref{ExamMnkRaBC} are justified by Theorem \ref{ThmLast} of Section \ref{secPF} which deals with a more extensive class of solutions to the quantum Yang--Baxter equation.

The notion of balanced oriented quantum algebra is analogous to the notion of quantum algebra as we shall see in the next section. To define regular isotopy invariants of oriented knots and links we will use an oriented quantum algebra with the additional structure of an invertible $G \in A$ which satisfies 
$$
t_{\sf d}(G) = t_{\sf u}(G) = G \qquad \mbox{and} \qquad t_{\sf d} {\circ} t_{\sf u}(x) = GxG^{-1}
$$ 
for all $x \in A$. The quintuple $(A, \rho, t_{\sf d}, t_{\sf u}, G)$ is a {\em twist oriented quantum algebra over} $k$, and the twist oriented quantum algebra $(A, \rho, 1_A, t_{\sf d} {\circ} t_{\sf u}, G)$ over $k$ is the {\em twist standard oriented quantum algebra associated with} $(A, \rho, t_{\sf d}, t_{\sf u}, G)$. When the underlying oriented quantum algebra structure of a twist oriented quantum algebra over $k$ is balanced we shall write $(A, \rho, t, G)$ for $(A, \rho, t_{\sf d}, t_{\sf u} , G)$, where $t = t_{\sf d} = t_{\sf u}$, and call $(A, \rho, t, G)$ a {\em twist balanced oriented quantum algebra over} $k$.

Balanced, or standard, quantum algebras have a twist structure in an important case. Let $A = {\rm M}_n(k)$ and suppose $t$ is an algebra automorphism  of $A$. By the Noether--Skolem Theorem there is an invertible $G \in A$ such that $t(x) = GxG^{-1}$ for all $x \in A$. See the corollary to
\cite[Theorem 4.3.1]{HER}. Observe that $G$ is unique up to scalar multiple since the center of $A$ is $k1$. 
\begin{LE}\label{LEMG}
Any balanced or any standard oriented quantum algebra structure $(A, \rho, t_{\sf d}, t_{\sf u})$ on $A = {\rm M}_n(k)$ extends to a twist oriented quantum algebra structure $(A, \rho, t_{\sf d}, t_{\sf u}, G)$ over $k$. Furthermore $G$ is unique up to scalar multiple. \qed
\end{LE}
\medskip

The reader can check that $({\rm M}_n(k), \rho_{a, {\sf B}, {\sf C}}, t, G)$ is a twist balanced oriented quantum algebra, where $G = \sum_{\imath = 1}^n \omega_\imath^2 E_{\imath \, \imath}$ and  $({\rm M}_n(k), \rho_{a, {\sf B}, {\sf C}}, t)$ is the balanced oriented quantum algebra of Example \ref{ExamMnkRaBC}. 
\section{Oriented Quantum Algebras and Quantum Algebras}\label{secOrientAlg}
Let $A$ be an algebra over $k$, let $\rho = \sum_{\imath = 1}^r a_\imath {\otimes} b_\imath \in A {\otimes} A$ be invertible with inverse described by $\rho^{-1} = \sum_{\jmath = 1}^s \alpha_\jmath {\otimes} \beta_\jmath$ and suppose that $t_{\sf d}, t_{\sf u}$ are algebra automorphisms of $A$. We begin this section by reformulating (qa.1)--(qa.3) in terms of these descriptions of $\rho$ and $\rho^{-1}$. 

We first observe that (qa.1) can be expressed 
\begin{equation}\label{Eqqaone}
\sum_{\imath = 1}^r \sum_{\jmath = 1}^s a_\imath t_{\sf d}(\alpha_\jmath ) {\otimes} \beta_\jmath t_{\sf u}(b_\imath ) = 1 {\otimes}1 = 
\sum_{\jmath = 1}^s \sum_{\imath = 1}^r t_{\sf d}(\alpha_\jmath )a_\imath  {\otimes} t_{\sf u}(b_\imath ) \beta_\jmath
\end{equation}
and (qa.2) is the same as 
\begin{equation}\label{Eqqatwo}
\sum_{\imath = 1}^r a_\imath {\otimes} b_\imath = \sum_{\imath = 1}^r t_{\sf d} (a_\imath) {\otimes} t_{\sf d}(b_\imath) = \sum_{\imath = 1}^r t_{\sf u} (a_\imath) {\otimes} t_{\sf u}(b_\imath). 
\end{equation}
Axiom (qa.3), which says that $\rho$ satisfies the quantum Yang--Baxter equation, is the familiar
\begin{equation}\label{Eqqathree}
\sum_{\imath, \jmath, \ell = 1}^r  a_\imath a_\jmath \otimes b_\imath a_\ell \otimes b_\jmath b_\ell =\sum_{\jmath, \imath, \ell = 1}^r  a_\jmath a_\imath \otimes a_\ell b_\imath \otimes b_\ell b_\jmath.
\end{equation}

Let $(A, \rho, t_{\sf d}, t_{\sf u})$ be an oriented quantum algebra. Using (\ref{Eqqaone}) and (\ref{Eqqathree}) it is not hard to see that $(A^{op},\rho, t_{\sf d}, t_{\sf u})$ is an oriented quantum algebra over $k$, which we denote by $A^{op}$ as well. Since $\rho$ is invertible and satisfies the quantum Yang--Baxter equation it follows that $\rho^{-1}$ does also. Let $t = t_{\sf d}$ or $t = t_{\sf u}$. Since $t{\otimes}t$ is an algebra automorphism of $A {\otimes} A$ and $\rho = (t{\otimes}t)(\rho)$ we have $\rho^{-1} = (t{\otimes}t)(\rho^{-1})$. Applying $t^{-1}_{\sf d} {\otimes} t^{-1}_{\sf u}$ to both sides of the equations of (\ref{Eqqaone}) we conclude that $(A, \rho^{-1}, t_{\sf d}^{-1}, t_{\sf u}^{-1})$ is an oriented quantum algebra over $k$. Observe that $(A, \rho^{op}, t_{\sf u}^{-1}, t_{\sf d}^{-1})$ is an oriented quantum algebra over $k$, where $\rho^{op} = \sum_{\imath = 1}^r b_\imath {\otimes} a_\imath$. If $K$ is a field extension of $k$ then $(A{\otimes}K, \rho{\otimes}1{\otimes}1, t_{\sf d}{\otimes}1_K, t_{\sf u}{\otimes}1_K)$ is an oriented quantum algebra over $K$, where we make the identification $\rho{\otimes}1{\otimes}1 = \sum_{\imath = 1}^r (a_\imath {\otimes}1){\otimes}(b_\imath {\otimes}1)$.

Suppose that $(A', \rho', t_{\sf d}', t_{\sf u}')$ is also an oriented quantum algebra over $k$ and write $\rho' = \sum_{\imath' = 1}^{r'} a'_{\imath'} {\otimes} b'_{\imath'} \in A' {\otimes} A'$. Then $(A {\otimes} A', \rho'', t_{\sf d} {\otimes} t_{\sf d}', t_{\sf u} {\otimes} t_{\sf u}')$ is an oriented quantum algebra over $k$, which we refer to as the {\em tensor product of} $(A, \rho, t_{\sf d}, t_{\sf u})$ {\em and} $(A', \rho', t_{\sf d}', t_{\sf u}')$, where $\rho'' = \sum_{\imath = 1}^r \sum_{\imath' = 1}^{r'} (a_\imath {\otimes} a'_{\imath'}) {\otimes} (b_\imath {\otimes} b'_{\imath'})$. Note that $(k, 1 {\otimes}1, 1_k)$ is a balanced oriented quantum algebra over $k$. 

A {\em morphism} $f : (A, \rho, t_{\sf d}, t_{\sf u}) \la (A', \rho', t_{\sf d}', t_{\sf u}')$ {\em of oriented quantum algebras} is an algebra map $ f : A \la A'$ which satisfies $\rho' = (f {\otimes} f)(\rho)$,  $t_{\sf d}' {\circ} f = f {\circ} t_{\sf d}$ and $t_{\sf u}' {\circ} f = f {\circ} t_{\sf u}$. Oriented quantum algebras together with their morphisms under composition form a monoidal category.

Suppose that $I$ is an ideal of $A$ and $t_{\sf d}(I) = t_{\sf u}(I) = I$. Then there is a unique oriented quantum algebra structure $(A/I, \overline{\rho}, \overline{t}_{\sf d}, \overline{t}_{\sf u})$ on the quotient algebra $A/I$ such that $\pi : (A, \rho, t_{\sf d}, t_{\sf u}) \la (A/I, \overline{\rho}, \overline{t}_{\sf d}, \overline{t}_{\sf u})$ is a morphism, where $\pi : A \la A/I$ is the projection. Furthermore, if $f : (A, \rho, t_{\sf d}, t_{\sf u}) \la (A', \rho', t_{\sf d}', t_{\sf u}')$ is a morphism of oriented quantum algebras and $f : A \la A'$ is onto, then $(A/{\rm ker}f, \overline{\rho}, \overline{t}_{\sf d}, \overline{t}_{\sf u})$ and $(A', \rho', t_{\sf d}', t_{\sf u}')$ are isomorphic oriented quantum algebras.

An {\em oriented quantum subalgebra of} $(A, \rho, t_{\sf d}, t_{\sf u})$ is an oriented quantum algebra $(B, \rho', t_{\sf d}', t_{\sf u}')$ over $k$, where $B$ is a subalgebra of $A$ and the inclusion $ \imath : B \la A$ determines a morphism $\imath : (B, \rho', t_{\sf d}', t_{\sf u}') \la (A, \rho, t_{\sf d}, t_{\sf u})$. In this case $\rho' = \rho$ and $t_{\sf d}(B) = t_{\sf u}(B) = B$. The inverse of $\rho'$ in $B{\otimes}B$ is necessarily $\rho^{-1}$; thus $\rho^{-1} \in B {\otimes} B$. Conversely, suppose that $B$ is a subalgebra of $A$ such that $\rho, \rho^{-1} \in B {\otimes} B$ and $t_{\sf d}(B) = t_{\sf u}(B) = B$. Then $(B, \rho, t_{\sf d}|_B, t_{\sf u}|_B)$ is an oriented quantum subalgebra of $(A, \rho, t_{\sf d}, t_{\sf u})$. 

The notion of minimal oriented quantum algebra is theoretically important in connection with regular isotopy invariants of oriented knots and links. A {\em minimal oriented quantum algebra} over $k$ is an oriented quantum algebra with exactly one oriented quantum subalgebra. We will show that any oriented quantum algebra has a unique minimal oriented quantum subalgebra.  Our proof is based on a bit of linear algebra.

Suppose that $V$ is a vector space over $k$, $\rho \in V {\otimes} V$ and let $V_{(\rho)} = \{ (u^* {\otimes}1_V)(\rho) + (1_V {\otimes} v^*)(\rho)\, | \, u^*, v^* \in V^*\}$. If $t$ is a linear endomorphism of $V$ which satisfies $\rho = (t {\otimes}t)(\rho)$ then $t(V_{(\rho)}) = V_{(\rho)}$.

To see this, we may assume that $\rho \neq 0$ and write $\rho = \sum_{\imath = 1}^r u_\imath {\otimes} v_\imath$, where $r$ is as small as possible. Then $\{ u_1, \ldots, u_r\}$, $\{ v_1, \ldots, v_r\}$ are linearly independent and the $u_\imath$'s together with the $v_\imath$'s span $V_{(\rho)}$. Since $\rho = (t {\otimes} t)(\rho)$, or equivalently $\sum_{\imath = 1}^r t(u_\imath) {\otimes} t(v_\imath) = \sum_{\imath = 1}^r u_\imath {\otimes} v_\imath$, it follows that the sets $\{ t(u_1), \ldots, t(u_r)\}$ and $\{ t(v_1), \ldots, t(v_r)\}$ are also linearly independent. It is easy to see now that $\{ u_1, \ldots, u_r\}$, $\{ t(u_1), \ldots, t(u_r)\}$ have the same span and that $\{ v_1, \ldots, v_r\}$, $\{ t(v_1), \ldots, t(v_r)\}$ have the same span. Thus $t(V_{(\rho)}) = V_{(\rho)}$.

By the definition of $V_{(\rho)}$ if $U$ is a subspace of $V$ such that $R \in U {\otimes} U$ then $V_{(\rho)} \subseteq U$. 
\begin{LE}\label{OrientedMin}
An oriented quantum algebra $(A, \rho, t_{\sf d}, t_{\sf u})$ over the field $k$ has a unique minimal oriented quantum subalgebra $(A_\rho, \rho, t_{\sf d}|_{A_\rho}, t_{\sf u}|_{A_\rho})$. Furthermore, if $(B, R, t_{\sf d}|_B, t_{\sf u}|_B)$ is any oriented quantum subalgebra of $(A, \rho, t_{\sf d}, t_{\sf u})$ then $A_\rho \subseteq B$.
\end{LE}

\pf 
Let $t = t_{\sf d}$ or $t = t_{\sf u}$. Since $\rho = (t {\otimes}t)(\rho)$ and $\rho^{-1} = (t {\otimes}t)(\rho^{-1})$ we conclude that $t(A_{(\rho)}) = A_{(\rho)}$ and $t(A_{(\rho^{-1})}) = A_{(\rho^{-1})}$ by the preceding remarks. Let $A_\rho$ be the subalgebra of $A$ generated by $A_{(\rho)} + A_{(\rho^{-1})}$. Then $t(A_\rho) = A_\rho$ and $\rho, \rho^{-1} \in A_\rho {\otimes} A_\rho$. Therefore $(A_\rho, \rho, t_{\sf d}|_{A_\rho}, t_{\sf u}|_{A_\rho})$ is an oriented quantum subalgebra of $(A, \rho, t_{\sf d}, t_{\sf u})$. If $(B, \rho, t_{\sf d}|_B, \rho, t_{\sf u}|_B)$ is any oriented quantum subalgebra of $(A, \rho, t_{\sf d}, t_{\sf u})$ then $\rho, \rho^{-1} \in B {\otimes} B$. Thus $A_{(\rho)}, A_{(\rho^{-1})} \subseteq B$, and hence $A_\rho \subseteq B$.
\qed
\medskip

Quantum algebras arise in the study of regular isotopy invariants of unoriented $1$--$1$ tangle,  knot and link diagrams; see 
\cite{HENN, KNOTS, RADKAUFFinv, RK2DIM, QC, KRS}.
A quantum algebra over $k$ is a triple $(A, \rho, s)$, where $A$ is an algebra over $k$, $R \in A \otimes A$ is invertible and $s : A \rTo A^{op}$ is an algebra isomorphism,  such that 
\begin{flushleft}
{\rm (QA.1)} 
$\rho^{-1} = (s \otimes 1_A)(\rho)$,
\vskip1\jot
{\rm (QA.2)} 
$\rho = (s \otimes s)(\rho)$ and
\vskip1\jot
{\rm (QA.3)}
$\rho_{1\, 2} \rho_{1\, 3} \rho_{2\, 3}  = \rho_{2\, 3} \rho_{1\, 3}\rho_{1\, 2}$.
\end{flushleft}
See
\cite[Section 3]{QC}
in particular. 

Let $(A, \rho, s)$ be a quantum algebra over $k$. Then $(A, \rho, s)$ has a unique minimal quantum subalgebra $(B, \rho, s|_B)$, and $B$ is generated by $A_{(\rho)} + A_{(\rho^{-1})} = A_{(\rho)}$ as an algebra. Suppose that $(A, \rho, t_{\sf d}, t_{\sf u})$ is an oriented quantum algebra over $k$.  Using (qa.2) we see that $t_{\sf d}(A_{(\rho)}) = t_{\sf u}(A_{(\rho)}) = A_{(\rho)}$. Therefore $t_{\sf d}(B) = t_{\sf u}(B) = B$ and $(B, \rho, t_{\sf d}|_B, t_{\sf u}|_B)$ is an oriented quantum subalgebra of $(A, \rho, t_{\sf d}, t_{\sf u})$. 

A good source of quantum algebras are the quasitriangular Hopf algebras. A quasitriangular Hopf algebra $(A, \rho)$ with antipode $s$ over the field $k$ has a quantum algebra structure $(A, \rho, s)$, and $(A, \rho, s)$ is minimal quantum algebra if and only if $(A, \rho)$ is a minimal quasitriangular Hopf algebra. See 
\cite[Section 2]{RMIN}.

Our next result characterizes all standard oriented quantum algebra structures on $A$ of the type $(A, \rho, 1_A, t)$ when $(A, \rho, s)$ is minimal.
\begin{PR}\label{MinQAoqa}
Let $(A, \rho, s)$ be a quantum algebra over the field $k$. Then: 
\begin{enumerate}
\item[{\rm a)}]
$(A, \rho, 1_A, s^{-2})$ is a standard oriented quantum algebra.
\item[{\rm b)}]
Suppose that $(A, \rho, s)$ is minimal and $(A, \rho, 1_A, t)$ is a standard oriented quantum algebra over $k$. Then $t = s^{-2}$.
\end{enumerate}
\end{PR}

\pf
Write $\rho = \sum_{\imath = 1}^r a_\imath {\otimes} b_\imath$. We first show part a). 

Since $(A, \rho, s)$ is a quantum algebra (qa.3) holds for $\rho$ and (qa.2) holds for $\rho$ and $s^{-2}$. 
The fact that $\rho^{-1} = \sum_{\imath = 1}^r s(a_{\imath}) {\otimes} b_\imath$ translates to 
$$
\sum_{\imath, \jmath = 1}^r s(a_\imath)a_{\jmath} {\otimes} b_\imath b_\jmath = 1 {\otimes} 1 = 
\sum_{\jmath, \imath = 1}^r a_{\jmath} s(a_\imath) {\otimes} b_\jmath b_\imath.
$$
Applying $s {\otimes} 1_A$ to both sides of these equations yields
$$
\sum_{\jmath, \imath = 1}^r s(a_{\jmath})s^2(a_\imath) {\otimes} b_\imath b_\jmath = 1 {\otimes} 1 = 
\sum_{\imath, \jmath = 1}^r s^2(a_\imath)s(a_{\jmath})  {\otimes} b_\jmath b_\imath.
$$
Since $\rho = (s^2 {\otimes} s^2)(\rho)$ it follows that 
$$
\sum_{\jmath, \imath = 1}^r s(a_{\jmath})a_\imath {\otimes} s^{-2}(b_\imath) b_\jmath = 1 {\otimes} 1 = 
\sum_{\imath, \jmath = 1}^r a_\imath s(a_{\jmath})  {\otimes} b_\jmath s^{-2}(b_\imath)
$$
which is to say that $(1_A {\otimes} s^{-2})(\rho)$ and $\rho^{-1}$ are inverses in $A{\otimes} A^{op}$. Therefore $(A, \rho, 1_A, s^{-2})$ is an oriented quantum algebra over $k$.

To show part b), suppose that $A = A_\rho$ and $(A, \rho, 1_A, t)$ is a standard oriented quantum algebra over $k$. Now $(1_A {\otimes} s^{-2})(\rho)$ is an inverse of $\rho^{-1}$ in $A {\otimes} A^{op}$ by part a). Thus  $(1_A {\otimes} s^{-2})(\rho) = (1_A {\otimes} t)(\rho)$. Since $t$ is one-one and $(t {\otimes} t)(\rho) = \rho = (s^{-2} {\otimes}s^{-2})(\rho)$ we have $(t {\otimes} 1_A)(\rho) = (s^{-2} {\otimes} 1_A)(\rho)$ also. Thus 
$$
\sum_{\imath = 1}^r a_\imath {\otimes} t(b_\imath) = 
\sum_{\imath = 1}^r a_\imath {\otimes} s^{-2}(b_\imath)
$$
and 
$$
\sum_{\imath = 1}^r t(a_\imath) {\otimes} b_\imath = 
\sum_{\imath = 1}^r s^{-2}(a_\imath) {\otimes} b_\imath.
$$
Assume that $r$ is as small as possible. Then $\{ a_1, \ldots, a_r\}$, $\{ b_1, \ldots, b_r\}$ are linearly independent. By virtue of the last two equations $t(b_\imath) = s^{-2}(b_\imath)$ and $t(a_\imath) = s^{-2}(a_\imath)$ for all $1 \leq \imath \leq r$. Thus $t$ and $s^{-2}$ agree on $A_{(\rho)}$ which generates $A_\rho = A$ as an algebra. Consequently $t = s^{-2}$ and part b) follows.
\qed
\medskip

By Proposition \ref{MinQAoqa} a minimal quantum algebra $(A, \rho, s)$ has a unique standard oriented quantum algebra structure of the form $(A, \rho, 1_A, t)$. Not every standard oriented quantum algebras over $k$ arises in this fashion by virtue of Proposition \ref{ToStandard} and the following example.
\begin{EX}
Let $n > 2$ and $({\rm M}_n(k), \rho_{a, {\sf B}, {\sf C}}, t)$ be the balanced oriented quantum algebra described by Example \ref{ExamMnkRaBC}. There is no quantum algebra of the form $({\rm M}_n(k), \rho_{a, {\sf B}, {\sf C}}, s)$.
\end{EX}
The assertion of the example follows by 
\cite[Proposition 1]{RKinv}.
By Proposition \ref{ToStandard} and part b) of Proposition \ref{MinQAoqa}:
\begin{CO}\label{tdtuEqualsminus2}
Suppose that $(A, \rho, s)$ is a minimal quantum algebra over $k$. If $(A, \rho, t_{\sf d}, t_{\sf u})$ is an oriented quantum algebra over $k$ then $t_{\sf d}{\circ} t_{\sf u} = s^{-2}$.  \qed
\end{CO}

Let $(A, \rho, s)$ be a minimal quantum algebra over $k$. Then $(A, \rho, s^{-2}, 1_A)$  and $(A, \rho, 1_A, s^{-2})$ are oriented quantum algebras over $k$ by Propositions \ref{ToStandard} and \ref{MinQAoqa}. There may be no other oriented quantum algebras of the type $(A, \rho, t)$. 

For example, suppose that the characteristic of $k$ is not $2$. Then Sweedler's $4$-dimensional Hopf algebra $A = A_{2, \, -1}$ over $k$ admits a family of minimal quasitriangular Hopf algebra structures $(A, \rho_\alpha )$ where $\alpha \in k^\star$. As an algebra $A$ is
generated over $k$ by symbols $a$ and $x$ which satisfy the relations
$$
a^2 = 1, \quad x^2 = 0 \quad \mbox{and} \quad xa = -ax.
$$
The coalgebra structure of $A$ is determined by
$$
\Delta (a) = a \otimes a \quad \mbox{and} \quad \Delta (x) = x \otimes a
+ 1 \otimes x.
$$
Consequently the antipode $s$ of $A$ satisfies $s(a) = a$ and $s(x) = ax$. Thus $s^2(x) = -x$ which means $s^2 \neq 1_A$. For each $\alpha \in k$ the element 
$$
\rho_{\alpha} = \frac{1}{2}(1 \otimes 1 + 1 \otimes a + a \otimes 1 - a
\otimes a) + \frac{\alpha}{2}(x \otimes x + x \otimes ax + ax \otimes
ax - ax \otimes x)
$$
gives $A$ the structure of a quasitriangular Hopf algebra
$(A, \rho_{\alpha})$ and every quasitriangular structure on $A$ is of this form. We remark that  $(A, \rho_{\alpha})$ is minimal quasitriangular if and only if $\alpha \neq 0$. For the details which justify these assertions the reader is referred to
\cite[Section 2]{RMIN}.
\begin{EX}
For $A = A_{2, \, -1}$ defined over a field $k$ of characteristic not $2$ the minimal quasitriangular Hopf algebra $(A, \rho_{\alpha})$ for all $\alpha \in k^\star$ admits only the oriented quantum algebra structures $(A, \rho_{\alpha}, 1_A, s^{-2})$ and $(A, \rho_{\alpha}, s^{-2}, 1_A)$.
\end{EX}

To see this, let $t$ be an algebra automorphism of $A$ which satisfies $\rho = (t{\otimes}t)(\rho)$. We will show that $t(a) = a$ and $t(x) = \pm x$. The assertion of the example will therefore follow by Corollary \ref{tdtuEqualsminus2}.  

It is easy to see that $t(a)$ is in the span of $a, x$, $ax$ and that $t(x)$ is in the span of $x$, $ax$. Thus $\rho = (t{\otimes}t)(\rho)$ if and only if $t(a) = a$ and $t(x) = \omega x + \rho ax$, where $\omega, \rho \in k$ and satisfy $\omega^2 + \rho^2 = 1$, $\omega^2 - \rho^2 + 2\omega \rho = 1$ and $\rho^2 - \omega^2 + 2\omega \rho = -1$. These equations hold if and only if  $\omega^2 = 1$ and $\rho = 0$. Therefore $f(a) = a$ and $f(x) = \pm x$. 

A quantum algebra $(A, \rho, s)$ over $k$ can be used to define a regular isotopy invariant of $1$--$1$ tangles. With the additional structure of an invertible $G \in A$ which satisfies $s(G) = G^{-1}$ and $s^2(x) = GxG^{-1}$ for all $x \in A$ a  regular isotopy invariant of knots and links can be defined. The quadruple $(A, \rho, s, G)$ is called a twist quantum algebra
\cite[Section 2]{RKinv}.
For the remainder of this section we shall be concerned with minimal quasitriangular Hopf algebras which admit a balanced oriented quantum algebra structure.
\begin{CO}\label{SMinusTwTtSq}
Let $(A, \rho)$ be a minimal quasitriangular Hopf algebra with antipode $s$ over the field $k$ and suppose that $t$ is a Hopf algebra automorphism of $A$. 
\begin{enumerate}
\item[{\rm a)}]
$(A, \rho, t)$ is a balanced oriented quantum algebra if and only if $\rho = (t{\otimes}t)(\rho)$ and $t^2 = s^{-2}$. 
\item[{\rm b)}]
Suppose that $(A, \rho, s, G)$ is a twist quantum algebra and $(A, \rho, t)$ is a balanced oriented quantum algebra. If $t(G) = G$ then $(A, \rho, t, G^{-1})$ is a twist balanced oriented quantum algebra.
\item[{\rm c)}]
Suppose that $(A, \rho, s^{2m})$ is a balanced oriented quantum algebra for some $m \geq 0$. Then $s^2$ has odd order.
\item[{\rm d)}]
Suppose that $s^2$ has order $2m + 1$ for some $m \geq 0$. Then $(A, \rho, s^{2m})$ is a balanced oriented quantum algebra.
\end{enumerate}
\end{CO}

\pf
Part b) follows from part a). Since $\rho = (s{\otimes}s)(\rho)$ and $s^2$ is a Hopf algebra automorphism of $A$, parts c) and d) follow from part a) also. Noting that $1_A {\otimes} t$ is an algebra automorphism of $A {\otimes} A^{op}$, part a) follows by Proposition \ref{MinQAoqa}.   
\qed
\medskip

Many of the finite-dimensional analogs of the quantized enveloping algebras have a Hopf algebra automorphism $t$ which satisfies $t^2 = s^{-2}$.

Now let $A$ be any finite-dimensional Hopf algebra with antipode $s$ over $k$ and consider the quantum double $D(A)$ defined in 
\cite{DRIN}.
As a coalgebra $D(A) = A^{*\, cop} {\otimes} A$. The multiplicative identity for the algebra structure on $D(A)$ is $\epsilon {\otimes}1$ and multiplication is determined by 
$$
(p {\otimes}a)(q{\otimes}b) = p(a_{(1)}{\rh}q{\lh}s^{-1}(a_{(3)})) {\otimes} a_{(2)}b
$$
for all $p, q \in A^*$ and $a, b \in A$; the functional $a{\rh}q{\lh}b  \in A^*$ is defined by $(a{\rh}q{\lh}b)(c) = p(bca)$ for all $c \in A$. Our description of the quantum double follows 
\cite{RMIN}.

Let $\{ a_1, \ldots, a_r\}$ be a linear basis for $A$ and let $\{ a^1, \ldots, a^r\}$ be the dual basis for $A^*$. Then $(D(A), \rho)$ is a minimal quasitriangular Hopf algebra, where $\rho = \sum_{\imath = 1}^r (\epsilon {\otimes} a_\imath ) {\otimes} (a^\imath {\otimes} 1)$. The definition of $\rho$ does not depend on the choice of basis for $A$. The square of antipode for $D(A)$ is $(s^*)^{-2} {\otimes} s^2$. We leave the proof of the following implication of part a) of  Corollary \ref{SMinusTwTtSq} to the reader.
\begin{CO}\label{DoubleOR}
Suppose that $A$ is a finite-dimensional Hopf algebra over the field $k$ and $t$ is a Hopf algebra automorphism of $A$. Then $(D(A), \rho, (t^{-1})^* {\otimes}t)$ is a balanced oriented quantum algebra if and only if  $t^2 = s^{-2}$. \qed
\end{CO}
\medskip

\noindent
Since the quantum double $(D(A), \rho)$ of a finite-dimensional Hopf algebra $A$ over $k$ is minimal quastriangular, $(D(A), \rho, 1_{D(A)}, (s^*)^2 {\otimes} s^{-2})$ is the unique standard oriented quantum algebra of the form $(D(A), \rho, 1_{D(A)}, t)$.

Taft's examples 
\cite{TAFT,TT},
which we denote by $A_{n, \, \omega}$, are generalizations of Sweedler's example and are made to order for producing Hopf algebra automorphisms $t$ which satisfy the condition of Corollary \ref{DoubleOR} when $k$ is an algebraically closed field of characteristic $0$. Generally $A_{n, \, \omega}$ is an $n^2$-dimensional Hopf algebra over $k$ and is defined as follows. Suppose $n \geq 2$ and $\omega \in k$ is a primitive $n^{th}$ root of unity. As a $k$-algebra $A_{n, \, \omega}$ generated by $a$ and $x$ subject to the relations
$$
a^n = 1, \quad x^n = 0 \quad \mbox{and} \quad xa = \omega ax.
$$
 The coalgebra structure of $A_{n, \, \omega}$ is determined by 
$$
\Delta (a) = a {\otimes} a \quad \mbox{and} \quad \Delta (x) = x {\otimes} a + 1 {\otimes} x
$$
and the antipode $s$ of $A_{n, \, \omega}$ is determined by $s(a) = a^{-1}$ and $s(x) = -xa^{-1}$. Thus $s^2(a) = a$ and $s^2(x) = axa^{-1} = \omega^{-1}x$. 

Now suppose $\beta \in k^\star$ is a square root of $\omega$ and let $t$ be the Hopf algebra automorphism of $A_{n, \, \omega}$ determined by $t(a) = a$ and $t(x) = \beta x$. Then $t^2(a) =  s^{-2}(a)$ and $t^2(x) = s^{-2}(x)$ which means $t^2 = s^{-2}$. A straightforward calculation shows that $(D(A_{n, \, \omega}), \rho, (t^{-1})^* {\otimes}t, \eta {\otimes} a^\ell)$ is a twist balanced oriented quantum algebra, where $0 \leq \ell < n$ and $\eta : A_{n, \, \omega} \la k$ is the algebra homomorphism determined by $\eta (a) = \omega^{-(\ell + 1)}$ and $\eta (x) = 0$.
\section{Invariants Associated with Oriented and Twist Oriented Quantum Algebras}\label{secInv}
We begin with an overview of the discussion of 
\cite{KRO}
regarding the regular isotopy invariants of oriented $1$--$1$ tangles determined by "bead sliding" and adaptation of this technique for the construction of regular isotopy invariants of oriented knots and links. We then develop a general theory of these invariants. In particular we show that these invariants can be computed from standard oriented quantum algebras.
\subsection{Invariants of Oriented $1$--$1$ Tangles Arising from Oriented Quantum Algebras}\label{secINV11T}
To describe the $1$--$1$ tangle invariants we begin with oriented diagrams. We represent oriented $1$--$1$ tangles as diagrams in the plane with respect to the vertical direction, for example 
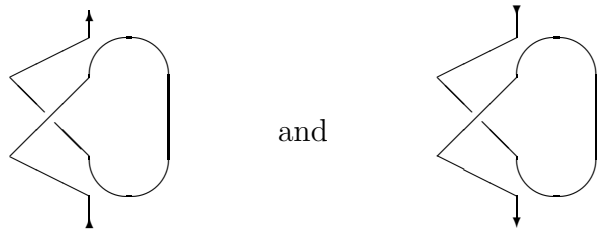
\begin{center}
\mbox{
\begin{picture}(80,90)(-10,-10)
\put(45,15){\oval(30,30)[b]}
\put(45,45){\oval(30,30)[t]}
\put(60,15){\line(0,1){30}}
\put(30,0){\line(0,-1){10}}
\put(30,-8){\vector(0,1){0}}
\put(30,60){\line(0,1){10}}
\put(30,70){\vector(0,1){0}}
%
%
\put(0,15){\line(1,1){30}}
\put(30,15){\line(-1,1){13}}
\put(0,45){\line(1,-1){13}}
%
%
\put(0,15){\line(2,-1){30}}
\put(0,45){\line(2,1){30}}
\end{picture}
}
\qquad 
\raisebox{6ex}{and}
\qquad 
\mbox{
\begin{picture}(80,90)(-10,-10)
\put(45,15){\oval(30,30)[b]}
\put(45,45){\oval(30,30)[t]}
\put(60,15){\line(0,1){30}}
\put(30,0){\line(0,-1){10}}
\put(30,-12){\vector(0,-1){0}}
\put(30,60){\line(0,1){10}}
\put(30,68){\vector(0,-1){0}}
%
%
\put(0,15){\line(1,1){30}}
\put(30,15){\line(-1,1){13}}
\put(0,45){\line(1,-1){13}}
%
%
\put(0,15){\line(2,-1){30}}
\put(0,45){\line(2,1){30}}
\end{picture}
}
\end{center}
\noindent
which we refer to as ${\bf T}_{\rm curl}$ and ${\bf T}^{op}_{\rm curl}$ respectively. The arrow heads indicate orientation.  Generally we require $1$--$1$ tangle diagrams to be completely contained in a box except for two protruding line segments as indicated in the two examples below.
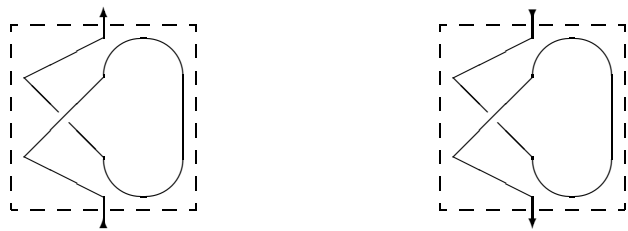
\begin{center}
\mbox{
\begin{picture}(80,90)(-10,-10)
\put(45,15){\oval(30,30)[b]}
\put(45,45){\oval(30,30)[t]}
\put(60,15){\line(0,1){30}}
\put(30,0){\line(0,-1){10}}
\put(30,-8){\vector(0,1){0}}
\put(30,60){\line(0,1){10}}
\put(30,72){\vector(0,1){0}}
%
%
\put(0,15){\line(1,1){30}}
\put(30,15){\line(-1,1){13}}
\put(0,45){\line(1,-1){13}}
%
%
\put(0,15){\line(2,-1){30}}
\put(0,45){\line(2,1){30}}
%
%
\put(-5,-5){\dashbox{5}(70,70)}
\end{picture}
}
\qquad \qquad \qquad 
\mbox{
\begin{picture}(80,90)(-10,-10)
\put(45,15){\oval(30,30)[b]}
\put(45,45){\oval(30,30)[t]}
\put(60,15){\line(0,1){30}}
\put(30,0){\line(0,-1){10}}
\put(30,-12){\vector(0,-1){0}}
\put(30,60){\line(0,1){10}}
\put(30,67){\vector(0,-1){0}}
%
%
\put(0,15){\line(1,1){30}}
\put(30,15){\line(-1,1){13}}
\put(0,45){\line(1,-1){13}}
%
%
\put(0,15){\line(2,-1){30}}
\put(0,45){\line(2,1){30}}
%
%
\put(-5,-5){\dashbox{5}(70,70)}
\end{picture}
}
\end{center}

Oriented $1$--$1$ tangle diagrams consist of some or all of the following components:
\begin{enumerate}
\item[$\bullet$]
oriented crossings;
\begin{enumerate}
\item[{\rm \phantom{a}}]
{\em under crossings}

%
%
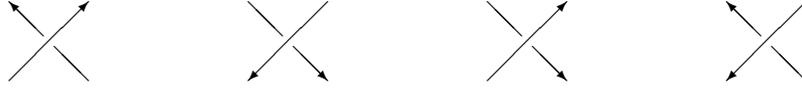
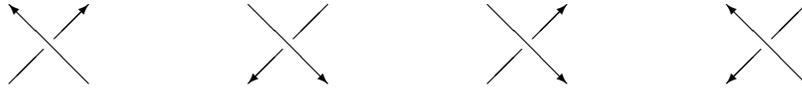
\mbox{
\begin{picture}(80,50)(-10,-10)
\put(0,0){\line(1,1){30}}
\put(0,30){\line(1,-1){13}}
\put(30,0){\line(-1,1){13}}
\put(0, 30){\vector(-1, 1){0}}
\put(30, 30){\vector(1, 1){0}}
\end{picture} }
%
%
\mbox{
\begin{picture}(80,50)(-10,-10)
\put(0,0){\line(1,1){30}}
\put(0,30){\line(1,-1){13}}
\put(30,0){\line(-1,1){13}}
\put(0, 0){\vector(-1, -1){0}}
\put(30, 0){\vector(1, -1){0}}
\end{picture} }
%
%
\mbox{
\begin{picture}(80,50)(-10,-10)
\put(0,0){\line(1,1){30}}
\put(0,30){\line(1,-1){13}}
\put(30,0){\line(-1,1){13}}
\put(30, 0){\vector(1, -1){0}}
\put(30, 30){\vector(1, 1){0}}
\end{picture} }
%
%
\mbox{
\begin{picture}(80,50)(-10,-10)
\put(0,0){\line(1,1){30}}
\put(0,30){\line(1,-1){13}}
\put(30,0){\line(-1,1){13}}
\put(0, 0){\vector(-1, -1){0}}
\put(0, 30){\vector(-1, 1){0}}
\end{picture} }
\item[{\rm \phantom{a}}]
{\em over crossings}

%
%
\mbox{
\begin{picture}(80,50)(-10,-10)
\put(0,0){\line(1,1){13}}
\put(30,0){\line(-1,1){30}}
\put(30,30){\line(-1,-1){13}}
\put(0, 30){\vector(-1, 1){0}}
\put(30, 30){\vector(1, 1){0}}
\end{picture} }
%
%
\mbox{
\begin{picture}(80,50)(-10,-10)
\put(0,0){\line(1,1){13}}
\put(30,0){\line(-1,1){30}}
\put(30,30){\line(-1,-1){13}}
\put(0, 0){\vector(-1, -1){0}}
\put(30, 0){\vector(1, -1){0}}
\end{picture} }
%
%
\mbox{
\begin{picture}(80,50)(-10,-10)
\put(0,0){\line(1,1){13}}
\put(0,30){\line(1,-1){30}}
\put(30,30){\line(-1,-1){13}}
\put(30, 0){\vector(1, -1){0}}
\put(30, 30){\vector(1, 1){0}}
\end{picture} }
%
%
\mbox{
\begin{picture}(80,50)(-10,-10)
\put(0,0){\line(1,1){13}}
\put(30,30){\line(-1,-1){13}}
\put(30,0){\line(-1,1){30}}
\put(0, 0){\vector(-1, -1){0}}
\put(0, 30){\vector(-1, 1){0}}
\end{picture} }
\end{enumerate}
\item[$\bullet$]
oriented local extrema;
\begin{enumerate}
\item[{\rm \phantom{a}}]
{\em local maxima}
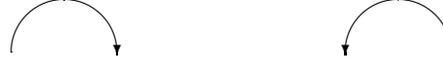
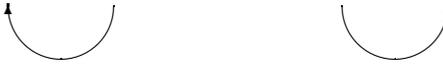
\mbox{
\begin{picture}(75,20)(-25,0)
\put(40,0){\vector(0, -1){}}
\put(20,0){\oval(40, 40)[t]}
\end{picture}}
\qquad 
\mbox{
\begin{picture}(75,20)(-45,0)
\put(0,0){\vector(0, -1){}}
\put(20,0){\oval(40, 40)[t]}
\end{picture}}

\item[{\rm \phantom{a}}]
{\em local minima}
\mbox{
\begin{picture}(75,40)(-25,-20)
\put(0,0){\vector(0, 1){}}
\put(20,0){\oval(40, 40)[b]}
\end{picture}}
\qquad  
\mbox{
\begin{picture}(75,40)(-45,-20)
\put(40,0){\vector(0, 1){}}
\put(20,0){\oval(40, 40)[b]}
\end{picture}}
\end{enumerate}
\end{enumerate}
and 
\begin{enumerate}
\item[$\bullet$]
oriented ``vertical" lines.
\end{enumerate}
The orientations of adjoining components of the tangle must be compatible.

If an oriented $1$--$1$ tangle diagram
%
%
\begin{picture}(16,26)(0,-5)
\put(8,-2){\vector(0, 1){0}}
\put(8,-5){\line(0,1){5}}
\put(8,24){\vector(0, 1){0}}
\put(8,16){\line(0,1){5}}
\put(4,4){${\bf T}$}
\put(0,0){\dashbox{4}(16,16)}
\end{picture}
(respectively
%
%
\begin{picture}(16,26)(0,-5)
\put(8,-7){\vector(0, -1){0}}
\put(8,-5){\line(0,1){5}}
\put(8,19){\vector(0, -1){0}}
\put(8,16){\line(0,1){5}}
\put(4,4){${\bf T}$}
\put(0,0){\dashbox{4}(16,16)}
\end{picture} $\;$)
can be decomposed into two oriented $1$--$1$ tangle diagrams
%
%
\begin{picture}(16,62)(0,-5)
\put(8,-2){\vector(0, 1){0}}
\put(8,-5){\line(0,1){5}}
\put(8,16){\line(0,1){10}}
\put(8,24){\vector(0, 1){0}}
\put(1,4){${\bf T}_1$}
\put(0,0){\dashbox{4}(16,16)}
\put(8,50){\vector(0, 1){0}}
\put(8,42){\line(0,1){5}}
\put(1,30){${\bf T}_2$}
\put(0,26){\dashbox{4}(16,16)}
\end{picture}
%
%
(respectively
\begin{picture}(16,62)(0,-5)
\put(8,-7){\vector(0, -1){0}}
\put(8,-5){\line(0,1){5}}
\put(8,16){\line(0,1){10}}
\put(8,19){\vector(0, -1){0}}
\put(1,4){${\bf T}_2$}
\put(0,0){\dashbox{4}(16,16)}
\put(8,46){\vector(0, -1){0}}
\put(8,42){\line(0,1){5}}
\put(1,30){${\bf T}_1$}
\put(0,26){\dashbox{4}(16,16)}
\end{picture}$\;$)
then we write ${\bf T} = {\bf T}_1 {\star} {\bf T}_2$.

Let ${\SL Tang}$ be the set of all oriented $1$--$1$ tangle diagrams with respect to the given vertical. For ${\bf T} \in {\SL Tang}$ we let ${\bf T}^{op}$ denote ${\bf T}$ with its orientation reversed.

Now suppose that $(A, \rho, t_{\sf d}, t_{\sf u})$ is an oriented quantum algebra over $k$ and let ${\bf T} \in {\SL Tang}$. We will construct an element $w({\bf T}) \in A$ which has the property $w({\bf T}) = w({\bf T}')$ whenever ${\bf T}, {\bf T}' \in {\SL Tang}$ are regularly isotopic. Thus the function ${\bf w}_A : {\SL Tang} \la A$ defined by 
$$
{\bf w}_A({\bf T}) = w({\bf T})
$$
for all ${\bf T} \in {\SL Tang}$ determines a regular isotopy invariant of oriented $1$--$1$ tangles.

To construct $w({\bf T})$ we will first construct a formal word $W({\bf T})$ as follows. If ${\bf T}$ has no crossings set $W({\bf T}) = 1$. Suppose that ${\bf T}$ has $n \geq 1$ crossings. Decorate the crossings of ${\bf T}$ according to the conventions of
\cite[Section 2.1]{KRO}, 
namely

%
%
\mbox{
\begin{picture}(80,50)(-10,-10)
\put(0,0){\line(1,1){30}}
\put(0,30){\line(1,-1){13}}
\put(30,0){\line(-1,1){13}}
\put(-5, 3){$e \; \bullet$}
\put(22, 3){$\bullet \; e'$}
\put(0, 30){\vector(-1, 1){0}}
\put(30, 30){\vector(1, 1){0}}
\end{picture} }
%
%
\mbox{
\begin{picture}(80,50)(-10,-10)
\put(0,0){\line(1,1){30}}
\put(0,30){\line(1,-1){13}}
\put(30,0){\line(-1,1){13}}
\put(23, 21){$\bullet \; e$}
\put(-9, 21){$e' \; \bullet$}
\put(0, 0){\vector(-1, -1){0}}
\put(30, 0){\vector(1, -1){0}}
\end{picture} }
%
%
\mbox{
\begin{picture}(80,50)(-10,-10)
\put(0,0){\line(1,1){30}}
\put(0,30){\line(1,-1){13}}
\put(30,0){\line(-1,1){13}}
\put(22, 21 ){$\bullet \; t_{\sf u}(E)$}
\put(22, 3){$\bullet \; E'$}
\put(30, 0){\vector(1, -1){0}}
\put(30, 30){\vector(1, 1){0}}
\end{picture} }
%
%
\mbox{
\begin{picture}(80,50)(-10,-10)
\put(0,0){\line(1,1){30}}
\put(0,30){\line(1,-1){13}}
\put(30,0){\line(-1,1){13}}
\put(-28, 3){$ t_{\sf d}(E) \; \bullet$}
\put(-13, 22){$E' \; \bullet$}
\put(0, 0){\vector(-1, -1){0}}
\put(0, 30){\vector(-1, 1){0}}
\end{picture} }

%
\mbox{
\begin{picture}(80,50)(-10,-10)
\put(0,0){\line(1,1){13}}
\put(30,0){\line(-1,1){30}}
\put(30,30){\line(-1,-1){13}}
\put(-9, 21){$E \; \bullet$}
\put(22, 21){$\bullet \; E'$}
\put(0, 30){\vector(-1, 1){0}}
\put(30, 30){\vector(1, 1){0}}
\end{picture} }
%
%
\mbox{
\begin{picture}(80,50)(-10,-10)
\put(0,0){\line(1,1){13}}
\put(30,0){\line(-1,1){30}}
\put(30,30){\line(-1,-1){13}}
\put(23, 3){$\bullet \; E$}
\put(-12, 3){$E' \; \bullet$}
\put(0, 0){\vector(-1, -1){0}}
\put(30, 0){\vector(1, -1){0}}
\end{picture} }
%
%
\mbox{
\begin{picture}(80,50)(-10,-10)
\put(0,0){\line(1,1){13}}
\put(0,30){\line(1,-1){30}}
\put(30,30){\line(-1,-1){13}}
\put(-5, 21){$e \; \bullet$}
\put(-27, 3){$t_{\sf u}(e') \; \bullet$}
\put(30, 0){\vector(1, -1){0}}
\put(30, 30){\vector(1, 1){0}}
\end{picture} }
%
%
\mbox{
\begin{picture}(80,50)(-10,-10)
\put(0,0){\line(1,1){13}}
\put(30,30){\line(-1,-1){13}}
\put(30,0){\line(-1,1){30}}
\put(22,3){$\bullet \; e$}
\put(22, 21){$\bullet \; t_{\sf d}(e')$}
\put(0, 0){\vector(-1, -1){0}}
\put(0, 30){\vector(-1, 1){0}}
\end{picture} }

\noindent
which one can view as diagrammatic representations of the formal expressions $\rho = e {\otimes} e'$, $\rho^{-1} =  E{\otimes} E'$, $(1_A {\otimes} t)(\rho) = e {\otimes} t(e')$ and $(t {\otimes} 1_A)(\rho^{-1}) = t(E) {\otimes} E'$, where $t = t_{\sf d}$ or $t = t_{\sf u}$. In practice we let $e {\otimes} e'$, $f {\otimes} f'$, $g {\otimes} g' \ldots$ denote copies of $\rho$ and $E {\otimes} E'$, $F {\otimes} F'$, $G {\otimes} G' \ldots$ denote copies of $\rho^{-1}$. Our labeled crossings are to be interpreted as flat diagrams of
\cite[Section 2.1]{KRO}
where the crossing type before decoration is indicated.

Think of the oriented tangle as a rigid wire and think of the decorations as labeled beads which slide freely around the wire. Starting at the beginning of the tangle diagram (in terms of orientation), traverse the diagram pushing the labeled beads along the wire so the that result is a juxtaposition of labeled beads at the end of the diagram. As a labeled bead passes through a local extremum its label is altered according to conventions of \cite[Section 2.1]{KRO}:
\begin{center}
\mbox{
\begin{picture}(75,40)(-25,-10)
\put(40,0){\vector(0, -1){}}
\put(20,0){\oval(40, 40)[t]}
\put(-10,5){$x \; \bullet$}
\end{picture}}
\qquad  
\raisebox{3ex}{to}
\qquad 
\mbox{
\begin{picture}(75,40)(-45,-10)
\put(40,0){\vector(0, -1){}}
\put(20,0){\oval(40, 40)[t]}
\put(35,5){$\bullet \; t_{\sf u}^{-1}(x)$}
\end{picture}}
\end{center}

\noindent
and
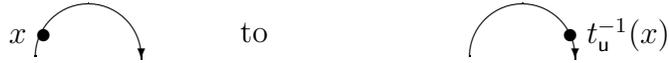
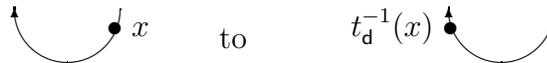
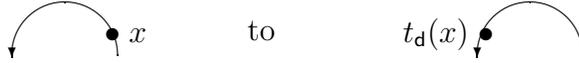
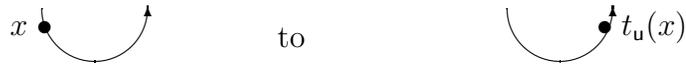
\begin{center}
\mbox{
\begin{picture}(75,40)(-25,-30)
\put(0,0){\vector(0, 1){}}
\put(20,0){\oval(40, 40)[b]}
\put(35,-10){$\bullet \; x$}
\end{picture}}
\qquad  
\raisebox{3ex}{to}
\qquad 
\mbox{
\begin{picture}(75,40)(-45,-30)
\put(0,0){\vector(0, 1){}}
\put(20,0){\oval(40, 40)[b]}
\put(-37,-10){$t_{\sf d}^{-1}(x) \; \bullet$}
\end{picture}}
\end{center}
for clockwise motion;
\begin{center}
\mbox{
\begin{picture}(75,40)(-25,-10)
\put(0,0){\vector(0, -1){}}
\put(20,0){\oval(40, 40)[t]}
\put(35,5){$\bullet \; x$}
\end{picture}}
\qquad  \quad
\raisebox{3ex}{to}
\qquad 
\mbox{
\begin{picture}(75,40)(-45,-10)
\put(0,0){\vector(0, -1){}}
\put(20,0){\oval(40, 40)[t]}
\put(-28,5){$t_{\sf d}(x) \;\bullet$}
\end{picture}}
\end{center}

\noindent
and
\begin{center}
\mbox{
\begin{picture}(75,40)(-25,-30)
\put(40,0){\vector(0, 1){}}
\put(20,0){\oval(40, 40)[b]}
\put(-12,-10){$x \; \bullet$}
\end{picture}}
\qquad  \quad
\raisebox{3ex}{to}
\qquad 
\mbox{
\begin{picture}(75,40)(-45,-30)
\put(40,0){\vector(0, 1){}}
\put(20,0){\oval(40, 40)[b]}
\put(34,-10){$\bullet \; t_{\sf u}(x)$}
\end{picture}}
\end{center}
for counterclockwise motion. We refer to the oriented local extrema
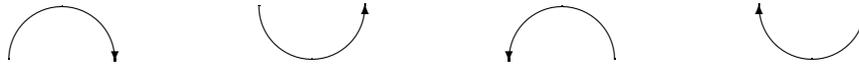
\begin{center}
\mbox{
\begin{picture}(75,40)(-25,-10)
\put(40,0){\vector(0, -1){}}
\put(20,0){\oval(40, 40)[t]}
\end{picture}}
\quad  
\mbox{
\begin{picture}(75,40)(-25,-30)
\put(40,0){\vector(0, 1){}}
\put(20,0){\oval(40, 40)[b]}
\end{picture}}
\quad
\mbox{
\begin{picture}(75,40)(-25,-10)
\put(0,0){\vector(0, -1){}}
\put(20,0){\oval(40, 40)[t]}
\end{picture}}
\quad
\mbox{
\begin{picture}(75,40)(-25,-30)
\put(0,0){\vector(0, 1){}}
\put(20,0){\oval(40, 40)[b]}
\end{picture}}
\end{center}
as having type (${\sf u}_-$), (${\sf u}_+$), (${\sf d}_+$) and  (${\sf d}_-$) respectively.
Reading the juxtaposed labeled beads in the direction of orientation results in a formal word $W({\bf T})$. Substituting copies of $\rho$ and $\rho^{-1}$ into their formal representations in $W({\bf T})$ results in an element $w({\bf T}) \in A$. 

The oriented $1$--$1$ tangle diagram ${\bf T}_{\rm trefoil}$ depicted below on the left is a good example to use for understanding the procedure for constructing  $w({\bf T})$.
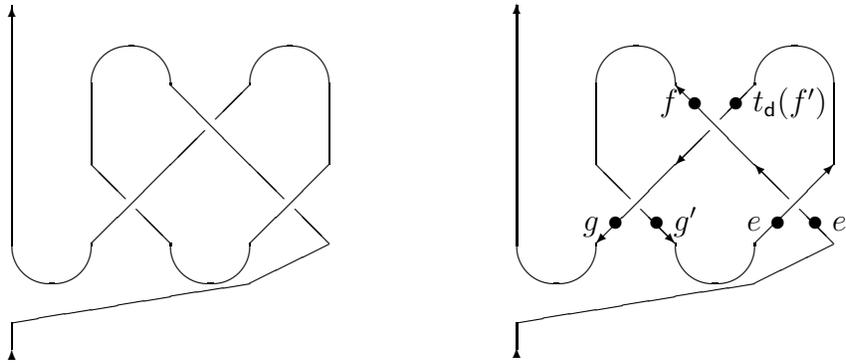
\begin{center}
\begin{picture}(140,140)(-10,-10)
\put(0,0){\line(0,-1){10}}
\put(0,-10){\vector(0, 1){0}}
\put(0,30){\line(0,1){90}}
\put(0,120){\vector(0, 1){0}}
\put(30,60){\line(0,1){30}}
\put(120,60){\line(0,1){30}}
\put(45,90){\oval(30,30)[t]}
\put(105,90){\oval(30,30)[t]}
\put(15,30){\oval(30,30)[b]}
\put(75,30){\oval(30,30)[b]}
\put(30,30){\line(1,1){30}}
\put(60,30){\line(-1,1){13}}
\put(30,60){\line(1,-1){13}}
\put(60,90){\line(1,-1){30}}
\put(60,60){\line(1,1){13}}
\put(90,90){\line(-1,-1){13}}
\put(90,30){\line(1,1){30}}
\put(90,60){\line(1,-1){13}}
\put(120,30){\line(-1,1){13}}
\put(120,30){\line(-2,-1){30}}
\put(0,0){\line(6,1){90}}
%
%
\end{picture}
\qquad \qquad
\begin{picture}(140,140)(-10,-10)
\put(0,0){\line(0,-1){10}}
\put(0,-9){\vector(0, 1){0}}
\put(0,30){\line(0,1){90}}
\put(0,121){\vector(0, 1){0}}
\put(30,60){\line(0,1){30}}
\put(120,60){\line(0,1){30}}
\put(45,90){\oval(30,30)[t]}
\put(105,90){\oval(30,30)[t]}
\put(15,30){\oval(30,30)[b]}
\put(75,30){\oval(30,30)[b]}
\put(30,30){\line(1,1){30}}
\put(25,35){$g \; \bullet$}
\put(60,30){\line(-1,1){13}}
\put(30,30){\vector(-1, -1){0}}
\put(30,60){\line(1,-1){13}}
\put(60,30){\vector(1, -1){0}}
\put(50,35){$\bullet \; g'$}
\put(60,90){\line(1,-1){30}}
\put(54,80){$f \; \bullet$}
\put(60,60){\vector(-1, -1){0}}
\put(60,60){\line(1,1){13}}
\put(90,90){\line(-1,-1){13}}
\put(60,90){\vector(-1, 1){0}}
\put(80,80){$\bullet \; t_{\sf d}(f')$}
\put(90,30){\line(1,1){30}}
\put(87,35){$e \; \bullet$}
\put(90,60){\vector(-1, 1){0}}
\put(90,60){\line(1,-1){13}}
\put(120,30){\line(-1,1){13}}
\put(120,60){\vector(1, 1){0}}
\put(110,35){$\bullet \; e'$}
\put(120,30){\line(-2,-1){30}}
\put(0,0){\line(6,1){90}}
%
%
\end{picture}
\end{center}
Traversal of the $1$--$1$ tangle diagram ${\bf T}_{\rm trefoil}$ results in the juxtaposition of labeled beads
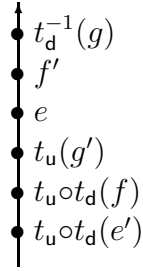
\begin{center}
\begin{picture}(40,110)(-10,-10)
\put(0,0){\line(0,-1){10}}
\put(0,0){\vector(0, 1){90}}
\put(-3.5,0){$\bullet \; t_{\sf u}{\circ} t_{\sf d}(e')$}
\put(-3.5,15){$\bullet \; t_{\sf u}{\circ} t_{\sf d}(f)$}
\put(-3.5,30){$\bullet \; t_{\sf u}(g')$}
\put(-3.5,45){$\bullet \; e$}
\put(-3.5,60){$\bullet \; f'$}
\put(-3.5,75){$\bullet \; t_{\sf d}^{-1}(g)$}
\end{picture}
\end{center}
Thus 
$$
W({\bf T}_{\rm trefoil}) = (t_{\sf u}{\circ} t_{\sf d})(e')(t_{\sf u}{\circ} t_{\sf d})(f)t_{\sf u}(g')ef't_{\sf d}^{-1}(g)
$$
from which we obtain after substitution 
$$
w({\bf T}_{\rm trefoil}) = \sum_{\imath, \jmath, \ell = 1}^r (t_{\sf u}{\circ} t_{\sf d})(b_\imath) (t_{\sf u}{\circ} t_{\sf d})(a_\jmath) t_{\sf u}(b_\ell)a_\imath b_\jmath t_{\sf d}^{-1}(a_\ell),
$$
where $\rho = \sum_{\imath = 1}^r a_\imath {\otimes} b_\imath \in A {\otimes} A$. Generally, the formal word $W({\bf T})$ can be viewed as merely a device which encodes instructions for defining an element 
$$
\rho^{\ell_1} {\otimes} \cdots {\otimes} \rho^{\ell_n} \in A {\otimes} \cdots {\otimes} A \;\; (2n \; \mbox{tensorands}),
$$
where $\ell_\imath = \pm 1$ for all $1 \leq \imath \leq n$, and describing the multilinear operations (permuting and applying powers of $t_{\sf d}, t_{\sf u}$ to tensorands) before the multiplication map  $A {\otimes} \cdots {\otimes} A \la A \;\; (a_1 {\otimes} \cdots {\otimes} a_{2n} \mapsto a_1\cdots a_{2n})$ is applied which results in $w({\bf T})$. 

Let $t = t_{\sf d}$ or $t = t_{\sf u}$. The axioms $\rho = (t {\otimes} t)(\rho)$ and $\rho^{-1} = (t {\otimes} t)(\rho^{-1})$, or symbolically $e {\otimes} e' = t(e) {\otimes} t(e')$ and $E {\otimes} E' = t(E) {\otimes} t(E')$ respectively, justify the rules  
$$
W({\bf T})  =  \cdots t^p(x) \cdots t^q(y) \cdots = \cdots t^{p+ \ell}(x) \cdots t^{q + \ell}(y) \cdots 
$$
for all integers $\ell$, where $x {\otimes} y$ or $y {\otimes} x$ represents either $\rho$ or $\rho^{-1}$. In light of these rules we may reformulate $w({\bf T}_{\rm trefoil})$ as 
$$
w({\bf T}_{\rm trefoil}) = \sum_{\imath, \jmath, \ell = 1}^r (t_{\sf u}{\circ} t_{\sf d}) (b_\imath) (t_{\sf u}{\circ} t_{\sf d}) (a_\jmath) (t_{\sf u}{\circ} t_{\sf d}) (b_\ell)a_\imath b_\jmath a_\ell.
$$

As a minor exercise the reader is encouraged to show that 
\begin{eqnarray*}
w({\bf T}^{op}_{\rm trefoil}) &  = & \sum_{\ell, \jmath, \imath  = 1}^r (t_{\sf u}^{-2}{\circ} t_{\sf d}^{-1}) (a_\ell) (t_{\sf u}^{-1}{\circ} t^{-1}_{\sf d}) (b_\jmath) (t_{\sf u}^{-1}{\circ} t^{-1}_{\sf d}) (a_\imath ) t^{-1}_{\sf u}(b_\ell) a_\jmath b_\imath \\
& = & \sum_{\ell, \jmath, \imath  = 1}^r a_\ell b_\jmath a_\imath (t_{\sf u}{\circ} t_{\sf d}) (b_\ell) (t_{\sf u}{\circ} t_{\sf d}) ( a_\jmath) (t_{\sf u}{\circ} t_{\sf d}) (b_\imath),
\end{eqnarray*}
and also that 
$$
w({\bf T}_{\rm curl}) = \sum_{\imath = 1}^r a_\imath (t_{\sf u}{\circ} t_{\sf d}) (b_\imath)  \quad \mbox{and} \quad 
w({\bf T}_{\rm curl}^{op}) = \sum_{\imath = 1}^r (t_{\sf u}{\circ} t_{\sf d})  (b_\imath) a_\imath.
$$

We now describe more precisely the procedure by which $W({\bf T})$ is calculated for the oriented  $1$--$1$ diagram ${\bf T}$ with $n \geq 1$ crossings. The formal word $W({\bf T})$ is the product of $2n$ factors, where each crossing contributes two factors in the manner described below.

Traverse the diagram ${\bf T}$ in the direction of orientation,  labeling the crossing lines of the diagram $1, \ldots, 2n$ in the order encountered. Denote by $u_{\sf d}(\ell)$ the number of local extrema which of type (${\sf d}_+$) minus the number of type (${\sf d}_-$) which are encountered in the part of the traversal of ${\bf T}$ from the line labeled $\ell$ to the end of the diagram. Define $u_{\sf u}(\ell)$ in the same manner, where (${\sf u}_+$) and (${\sf u}_-$) replace (${\sf d}_+$) and (${\sf d}_-$) respectively. 

Let $\chi$ be a crossing and suppose that its lines are labeled $\imath$ and $\jmath$, where $\imath < \jmath$, and let $x$ and $y$ be the decorations on the crossing lines $\imath$ and $\jmath$ respectively. Then the contribution which $\chi$ makes to $W({\bf T})$ is 
$$
W({\bf T})   =   \cdots t_{\sf d}^{u_{\sf d}(\imath)}{\circ} t_{\sf u}^{u_{\sf u}(\imath)}(x) \cdots  t_{\sf d}^{u_{\sf d}(\jmath)}{\circ} t_{\sf u}^{u_{\sf u}(\jmath)}(y) \cdots  \\
$$
where the indicated factors are the $\imath^{th}$ and $\jmath^{th}$ respectively of the product.

Let ${\bf T}, {\bf T}' \in {\SL Tang}$. To show that ${\bf w}_A({\bf T}) = {\bf w}_A({\bf T}')$ whenever ${\bf T}, {\bf T}'$ are regularly isotopic is a matter of showing that ${\bf w}_A({\bf T})$ is unaffected when a local portion of ${\bf T}$ is replaced by its equivalent according to 
%
\begin{flushleft}
\raisebox{6ex}{(M.1) \quad}
\mbox{
\begin{picture}(60,50)(-30,-25)
\put(-30,0){\line(6,-5){30}}
\put(0,-25){\line(0,-1){5}}
\put(0,25){\line(0,1){5}}
\put(30,0){\line(-6,5){30}}
\put(-15, 0){\oval(30,30)[t]}
\put(15, 0){\oval(30,30)[b]}
\end{picture} }
\raisebox{4ex}{$\;\; \approx \;\;$}
\mbox{
\begin{picture}(10,60)(-5,0)
\put(0,-5){\line(0,1){60}}
\end{picture} }
\raisebox{4ex}{\quad and \quad}
%
%
\mbox{
\begin{picture}(60,50)(-30,-25)
\put(-30,0){\line(6,5){30}}
\put(0,-25){\line(0,-1){5}}
\put(0,25){\line(0,1){5}}
\put(30,0){\line(-6,-5){30}}
\put(-15, 0){\oval(30,30)[b]}
\put(15, 0){\oval(30,30)[t]}
\end{picture} }
\raisebox{4ex}{$\;\; \approx \;\;$}
\mbox{
\begin{picture}(10,60)(-5,0)
\put(0,-5){\line(0,1){60}}
\end{picture} }
\end{flushleft}
\vskip1\jot
%
\begin{flushleft}
\raisebox{6ex}{(M.2) \quad}
\mbox{
\begin{picture}(30,70)(0,-5)
\put(0,0){\line(0,-1){5}}
\put(0,60){\line(0,1){5}}
\put(30,0){\line(0,-1){5}}
\put(30,60){\line(0,1){5}}
\put(0,0){\line(1,1){30}}
\put(0,30){\line(1,-1){13}}
\put(30,0){\line(-1,1){13}}
\put(0,30){\line(1,1){13}}
\put(30,60){\line(-1,-1){13}}
\put(0,60){\line(1,-1){30}}
\end{picture} }
\raisebox{6ex}{$\;\; \approx \;\; $}
\mbox{
\begin{picture}(30,70)(0,-5)
\put(0,-5){\line(0,1){70}}
\put(30,-5){\line(0,1){70}}
\end{picture} }
\end{flushleft}
\vskip1\jot
\begin{flushleft}
\raisebox{9ex}{(M.3) \quad}
\mbox{
\begin{picture}(60,100)(0,-5)
\put(0,0){\line(0,-1){5}}
\put(30,0){\line(0,-1){5}}
\put(60,0){\line(0,-1){5}}
\put(0,90){\line(0,1){5}}
\put(30,90){\line(0,1){5}}
\put(60,90){\line(0,1){5}}
\put(0,30){\line(0,1){30}}
\put(60,0){\line(0,1){30}}
\put(60,60){\line(0,1){30}}
\put(0,0){\line(1,1){13}}
\put(30,30){\line(-1,-1){13}}
\put(0,30){\line(1,-1){30}}
\put(0,60){\line(1,1){13}}
\put(30,90){\line(-1,-1){13}}
\put(0,90){\line(1,-1){30}}
\put(30,30){\line(1,1){13}}
\put(60,60){\line(-1,-1){13}}
\put(30,60){\line(1,-1){30}}
\end{picture} }
\raisebox{9ex}{$\;\; \approx \;\;$}
\mbox{
\begin{picture}(60,100)(0,-5)
\put(0,0){\line(0,-1){5}}
\put(30,0){\line(0,-1){5}}
\put(60,0){\line(0,-1){5}}
\put(0,90){\line(0,1){5}}
\put(30,90){\line(0,1){5}}
\put(60,90){\line(0,1){5}}
\put(0,0){\line(0,1){30}}
\put(0,60){\line(0,1){30}}
\put(60,30){\line(0,1){30}}
\put(0,30){\line(1,1){13}}
\put(30,60){\line(-1,-1){13}}
\put(30,30){\line(-1,1){30}}
\put(30,0){\line(1,1){13}}
\put(60,30){\line(-1,-1){13}}
\put(60,0){\line(-1,1){30}}
\put(30,60){\line(1,1){13}}
\put(60,90){\line(-1,-1){13}}
\put(60,60){\line(-1,1){30}}
\end{picture} }
\end{flushleft}
\vskip1\jot
%
\begin{flushleft}
\raisebox{9ex}{(M.4) \quad}
\mbox{
\begin{picture}(100,100)(0,0)
\put(0,30){\line(0,-1){5}}
\put(45,15){\line(0,-1){5}}
\put(90,30){\line(0,-1){5}}
\put(45,90){\line(0,1){5}}
\put(0,60){\line(3,2){45}}
\put(60,60){\line(1,-1){30}}
\put(30,30){\line(1,-1){15}}
\put(0,30){\line(1,1){30}}
\put(0,60){\line(1,-1){13}}
\put(30,30){\line(-1,1){13}}
\put(45,60){\oval(30,30)[t]}
\end{picture} }
\raisebox{9ex}{$\;\; \approx \;\;$}
\mbox{
\begin{picture}(100,100)(0,0)
\put(0,30){\line(0,-1){5}}
\put(45,15){\line(0,-1){5}}
\put(90,30){\line(0,-1){5}}
\put(45,90){\line(0,1){5}}
\put(90,60){\line(-3,2){45}}
\put(0,30){\line(1,1){30}}
\put(60,30){\line(-1,-1){15}}
\put(60,60){\line(1,-1){30}}
\put(60,30){\line(1,1){13}}
\put(90,60){\line(-1,-1){13}}
\put(45,60){\oval(30,30)[t]}
\end{picture} }
\end{flushleft}

and
%
%
\begin{flushleft}
\raisebox{9ex}{\phantom{(M.4)} \quad}
\mbox{
\begin{picture}(100,100)(0,0)
\put(0,60){\line(0,1){5}}
\put(45,75){\line(0,1){5}}
\put(90,60){\line(0,1){5}}
\put(45,0){\line(0,-1){5}}
\put(0,30){\line(3,-2){45}}
\put(60,30){\line(1,1){30}}
\put(30,60){\line(1,1){15}}
\put(0,60){\line(1,-1){30}}
\put(0,30){\line(1,1){13}}
\put(30,60){\line(-1,-1){13}}
\put(45,30){\oval(30,30)[b]}
\end{picture} }
\raisebox{6ex}{$\;\; \approx \;\;$}
\mbox{
\begin{picture}(100,100)(0,0)
\put(0,60){\line(0,1){5}}
\put(45,75){\line(0,1){5}}
\put(90,60){\line(0,1){5}}
\put(45,0){\line(0,-1){5}}
\put(90,30){\line(-3,-2){45}}
\put(60,60){\line(-1,1){15}}
\put(0,60){\line(1,-1){30}}
\put(60,30){\line(1,1){30}}
\put(60,60){\line(1,-1){13}}
\put(90,30){\line(-1,1){13}}
\put(45,30){\oval(30,30)[b]}
\end{picture} }
\end{flushleft}
\vskip1\jot

\noindent
and (M.2rev)--(M.4rev), which are (M.2)--(M.4) respectively with overcrossing lines replaced by under crossing lines and vice versa. Observe that the "twist" equivalences described below are consequences of (M.1) and (M.4).
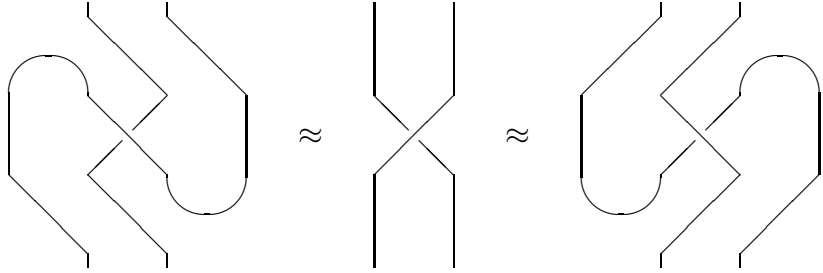
\begin{center}
%
\mbox{
\begin{picture}(100,100)(-5,-5)
\put(30,0){\line(0,-1){5}}
\put(60,0){\line(0,-1){5}}
\put(30,90){\line(0,1){5}}
\put(60,90){\line(0,1){5}}
\put(0,30){\line(1,-1){30}}
\put(30,30){\line(1,-1){30}}
\put(30,90){\line(1,-1){30}}
\put(60,90){\line(1,-1){30}}
\put(15,60){\oval(30,30)[t]}
\put(75,30){\oval(30,30)[b]}
\put(0,30){\line(0,1){30}}
\put(90,30){\line(0,1){30}}
\put(30,30){\line(1,1){13}}
\put(30,60){\line(1,-1){30}}
\put(60,60){\line(-1,-1){13}}
\end{picture} }
\raisebox{9ex}{$\;\; \approx \;\;$}
\mbox{
\begin{picture}(40,100)(-5,-5)
\put(0,0){\line(0,-1){5}}
\put(30,0){\line(0,-1){5}}
\put(0,90){\line(0,1){5}}
\put(30,90){\line(0,1){5}}
\put(0,0){\line(0,1){30}}
\put(0,60){\line(0,1){30}}
\put(30,0){\line(0,1){30}}
\put(30,60){\line(0,1){30}}
\put(0,30){\line(1,1){30}}
\put(0,60){\line(1,-1){13}}
\put(30,30){\line(-1,1){13}}
\end{picture} }
%
%
\raisebox{9ex}{$\;\; \approx \;\;$}
%
\mbox{
\begin{picture}(100,100)(-5,-5)
\put(30,0){\line(0,-1){5}}
\put(60,0){\line(0,-1){5}}
\put(30,90){\line(0,1){5}}
\put(60,90){\line(0,1){5}}
\put(30,0){\line(1,1){30}}
\put(60,0){\line(1,1){30}}
\put(0,60){\line(1,1){30}}
\put(30,60){\line(1,1){30}}
\put(15,30){\oval(30,30)[b]}
\put(75,60){\oval(30,30)[t]}
\put(0,30){\line(0,1){30}}
\put(90,30){\line(0,1){30}}
\put(30,30){\line(1,1){13}}
\put(30,60){\line(1,-1){30}}
\put(60,60){\line(-1,-1){13}}
\end{picture} }
\end{center}
These twist equivalences, and their counterparts derived from (M.1) and (M.4rev), are sometimes useful for the calculation of ${\bf w}_A({\bf T})$ in that we may assume all crossing lines are pointed in the same direction -- up, down, to the right or to the left.

Let $1 \leq \ell \leq 2n$ and let $d(\ell)$ be the Whitney degree of the portion of the diagram ${\bf T}$ traversed from line $\ell$ to the end. If this traversal starts and ends in the same direction then $u_{\sf d}(\ell) = u_{\sf u}(\ell) = -d(\ell)$. For denote the number of local extrema of the type (${\sf d}_+$), (${\sf d}_-$), (${\sf u}_+$) and $({\sf u}_-$) encountered on the traversal from line $\ell$ to the end of ${\bf T}$ by ${\sf d}_+({\bf T} : \ell), {\sf d}_-({\bf T} : \ell), {\sf u}_+({\bf T} : \ell)$ and ${\sf u}_+({\bf T} : \ell)$ respectively. Then 
$$
{\sf d}_-({\bf T} : \ell) + {\sf u}_-({\bf T} : \ell) - {\sf d}_+({\bf T} : \ell) - {\sf u}_+({\bf T} : \ell)
= 2d(\ell).
$$
Since the traversal of ${\bf T}$ from line $\ell$ starts and ends in the same direction it follows that $
{\sf d}_+({\bf T} : \ell) + {\sf u}_-({\bf T} : \ell) = {\sf d}_-({\bf T} : \ell) + {\sf u}_+({\bf T} : \ell)$. Consequently 
$$
d(\ell) = {\sf d}_-({\bf T} : \ell) - {\sf d}_+({\bf T} : \ell) = {\sf u}_-({\bf T} : \ell) - {\sf u}_+({\bf T} : \ell)
$$ 
and thus $u_{\sf d}(\ell) = u_{\sf u}(\ell) = -d(\ell)$.

The reader may have noticed that $w({\bf T}^{op})$ is $w({\bf T})$ with its factors reversed when ${\bf T} = {\bf T}_{\rm trefoil}$ and  ${\bf T} = {\bf T}_{\rm curl}$. 
\begin{PR}\label{wTop}
Let $A$ be an oriented quantum algebra over the field $k$. Then ${\bf w}_A({\bf T}^{op}) = {\bf w}_{A^{op}}({\bf T})$ for all ${\bf T} \in {\SL Tang}$
\end{PR}

\pf 
We assume that ${\bf T}$ has $n \geq 1$ crossings. Let $\chi$ be a crossing of ${\bf T}$, let $\imath$ and $\jmath$, where $\imath < \jmath$, be the labels for the lines of $\chi$ which result from a traversal of ${\bf T}$ in the direction of orientation and let $x$ and $y$ be the decorations on the crossing lines labeled $\imath$ and $\jmath$ respectively. Then $\chi$ contributes factors $t_{\sf d}^{u_{\sf d}(\imath)} {\circ}t_{\sf u}^{u_{\sf u}(\imath)}(x)$ and $t_{\sf d}^{u_{\sf d}(\jmath)} {\circ}t_{\sf u}^{u_{\sf u}(\jmath)}(y)$ to $W({\bf T})$ which are located in positions $\imath$ and $\jmath$ respectively. 

Observe that the lines labeled $\imath$ and $\jmath$ in the traversal of ${\bf T}$ are labeled $2n + 1 - \imath$ and $2n + 1 - \jmath$ respectively in the traversal of ${\bf T}^{op}$. Let $s_{\sf d}$ be the number of  local extrema of ${\bf T}$ of type (${\sf d}_+$) minus the number of type (${\sf d}_-$) and let $s_{\sf u}$ be the number of type (${\sf u}_+$) minus the number of type (${\sf u}_-$). Then $u_{\sf d}^{op}(2n + 1 - \ell) = u_{\sf u}(\ell) - s_{\sf u}$ and $u_{\sf u}^{op}(2n + 1 - \ell) =  u_{\sf d}(\ell) - s_{\sf d}$ for all $1 \leq \ell \leq 2n$. Thus
\begin{eqnarray*}
W({\bf T}^{op}) & = & \cdots t_{\sf d}^{u_{\sf d}^{op}(2n + 1 - \jmath)} {\circ} t_{\sf u}^{u_{\sf u}^{op}(2n + 1 - \jmath)}(y) \cdots t_{\sf d}^{u_{\sf d}^{op}(2n + 1 - \imath)} {\circ} t_{\sf u}^{u_{\sf u}^{op}(2n + 1 - \imath)}(x) \cdots \\
& = & \cdots t_{\sf d}^{u_{\sf u}(\jmath) - s_{\sf u}} {\circ} t_{\sf u}^{u_{\sf d}(\jmath) - s_{\sf d}}(y) \cdots t_{\sf d}^{u_{\sf u}(\imath) - s_{\sf u}} {\circ} t_{\sf u}^{u_{\sf d}(\imath) - s_{\sf d}}(x) \cdots \\
& = &\cdots t_{\sf d}^{u_{\sf u}(\jmath)} {\circ} t_{\sf u}^{u_{\sf d}(\jmath)}(y) \cdots t_{\sf d}^{u_{\sf u}(\imath) } {\circ} t_{\sf u}^{u_{\sf d}(\imath) }(x) \cdots
\end{eqnarray*}
where the factors are in positions $2n + 1 - \jmath$ and $2n + 1 - \imath$ respectively. Now to prove the proposition we may assume that all crossings of ${\bf T}$ are directed upward. In this case $u_{\sf d}(\ell) = u_{\sf u}(\ell)$ for all $1 \leq \ell \leq 2n$ and which establishes ${\bf w}_A({\bf T}^{op}) = {\bf w}_{A^{op}}({\bf T})$.
\qed
\begin{PR}
Suppose that $f : A \la B$ is a morphism of oriented quantum algebras. Then $f({\bf w}_A({\bf T})) = {\bf w}_B({\bf T})$ for all ${\bf T} \in {\SL Tang}$. 
\end{PR}
\qed
\medskip

The proposition, the proof of which is easy, makes an interesting statement about regular isotopy invariants of the type ${\bf w}_A$. Whenever $A$ and $B$ are oriented quantum algebras related by a morphism $f : A \la B$ then ${\bf w}_A$ dominates ${\bf w}_B$; that is if ${\bf T},  {\bf T}' \in {\SL Tang}$ satisfy ${\bf w}_A({\bf T}) = {\bf w}_A({\bf T}')$ then ${\bf w}_B({\bf T}) = {\bf w}_B({\bf T}')$.

We shall also write ${\bf Inv}_A$ for ${\bf w}_A$ to be consistent with the notation of the next section. If ${\bf T} \in {\SL Tang}$ can be written ${\bf T} = {\bf T}_1{\star} {\bf T}_2$, note that ${\bf w}_A({\bf T}) = {\bf w}_A({\bf T}_1) {\bf w}_A({\bf T}_2)$ or equivalently ${\bf Inv}_A({\bf T}_1{\star} {\bf T}_2) = {\bf Inv}_A({\bf T}_1) {\bf Inv}_A({\bf T}_2)$.

By virtue of the next result the $1$--$1$ tangle invariants described in this section are computed by standard oriented quantum algebras.
\begin{TH}\label{11Standard}
Let $A$ be an oriented quantum algebra over the field $k$ and suppose that $A_s$ is the associated standard oriented quantum algebra. Then ${\bf Inv}_A({\bf T}) = {\bf Inv}_{A_s}({\bf T})$ for all ${\bf T} \in {\SL Tang}$
\end{TH}

\pf
Let ${\bf T} \in {\SL Tang}$. To show that ${\bf Inv}_A({\bf T}) = {\bf Inv}_{A_s}({\bf T})$ we may assume that ${\bf T}$ has at least one crossing and that all crossings of ${\bf T}$ are directed upward. Let $\chi$ be a crossing of ${\bf T}$ with line labels $\imath$ and $\jmath$, where $\imath < \jmath$. Then 
\begin{eqnarray*}
W({\bf T}) &  = & \cdots t_{\sf d}^{u_{\sf d}(\imath)} {\circ} t_{\sf u}^{u_{\sf u}(\imath)}(x) \cdots 
\cdots t_{\sf d}^{u_{\sf d}(\jmath)} {\circ} t_{\sf u}^{u_{\sf u}(\jmath)}(y) \cdots \\
& = & \cdots (t_{\sf d} {\circ} t_{\sf u})^{-d(\imath)} (x) \cdots (t_{\sf d} {\circ} t_{\sf u})^{-d(\jmath)} (y) \cdots \; .
\end{eqnarray*}
\qed
\subsection{Invariants of Oriented Knots and Links Arising from Twist Oriented Quantum Algebras}\label{KnotLinkINV}
Let ${\SL Knot}$ (respectively ${\SL Link}$) denote the set of oriented knot (respectively link) diagrams with respect to a fixed vertical in the plane. Throughout this section $(A, \rho, t_{\sf d}, t_{\sf u}, G)$ is a twist oriented quantum algebra over the field $k$. Notice that $(A^{op}, R, t_{\sf d}, t_{\sf u}, G^{-1})$ is an oriented quantum algebra over $k$, which we denote by $A^{op}$.  

If $B$ is any algebra over $k$ then  ${\sf tr} \in B^*$ is a {\em tracelike element } if ${\sf tr}(ab) = {\sf tr}(ba)$ for all $a, b \in B$. For any tracelike element ${\sf tr} \in A^*$ which is $t_{\sf d}^*, t_{\sf u}^*$-invariant, that is satisfies ${\sf tr} {\circ}t_{\sf d} = {\sf tr} {\circ}t_{\sf u}  = {\sf tr}$,  we construct a function ${\bf Inv}_{A, \, {\sf tr}} :  {\SL Link} \la k$ with the property that ${\bf L}, {\bf L}' \in {\SL Link}$ regularly isotopic implies ${\bf Inv}_{A, \, {\sf tr}}({\bf L}) = {\bf Inv}_{A, \, {\sf tr}}({\bf L}')$. Thus ${\bf Inv}_{A, \, {\sf tr}}$ determines a regular isotopy invariant of oriented links. To begin we define ${\bf Inv}_{A, \, {\sf tr}}({\bf K})$ for ${\bf K} \in {\SL Knot}$.

Let ${\bf K} \in {\SL Knot}$. To define ${\bf Inv}_{A, \, {\sf tr}}({\bf K})$  we first construct an element $w({\bf K}) \in A$. If ${\bf K}$ has no crossings we set $w({\bf K}) = 1$. 

Suppose that ${\bf K}$ has $n \geq 1$ crossings. Decorate the crossings of ${\bf K}$ according to the conventions of Section \ref{secINV11T} and choose a point $P$ on a vertical line in the knot diagram ${\bf K}$. (There is no harm, under regular isotopy, in inserting a vertical line at the end of a crossing line or local extrema -- thus we may assume that ${\bf K}$ has a vertical line.) 

Traverse the knot diagram ${\bf K}$, starting at $P$ and moving in the direction of the orientation back to $P$, and label the crossing lines $1, \ldots, 2n$ in the order encountered. Let $W({\bf K})$ be a formal product with $2n$ factors, where each crossing contributes two factors according to the algorithm described in Section \ref{secINV11T} for oriented $1$--$1$ tangle diagrams, and let $w({\bf K})$ be obtained from $W({\bf K})$ in the same manner that $w({\bf T})$ was obtained from $W({\bf T})$ in Section \ref{secINV11T}, where ${\bf T} \in {\SL Tang}$. Here ``to the end of the (tangle) diagram" is replaced by ``back to $P$ in the direction of orientation".

Let $d$ be the Whitney degree of the oriented knot diagram ${\bf K}$. Then $2d$ is the number of local extrema with clockwise orientation minus the number of extrema with counterclockwise orientation. We will show that the scalar ${\sf tr}(G^d w({\bf K}))$ does not depend on the starting point $P$. 

Consider a new starting point $P_{new}$ which follows $P$ in the orientation of ${\bf K}$.  Let $W_{new}({\bf K})$ and $w_{new}({\bf K})$ for $P_{new}$ be the analogs of $W({\bf K})$ and $w({\bf K})$ respectively for $P$. There are two cases to consider. 

Suppose first of all that traversal of the diagram ${\bf K}$ from $P$ to $P_{new}$ in the direction of orientation passes through exactly one local extremum and no crossing lines. By examining the four local extremum types one sees that there are $r_{\sf d}, r_{\sf u} \in \{ -1, 0, 1\}$ such that $u_{{\sf d}, new}(\ell) = u_{\sf d}(\ell) + r_{\sf d}$ and $u_{{\sf u}, new}(\ell) = u_{\sf u}(\ell) + r_{\sf u}$ for all $1 \leq \ell \leq 2n$. Since $\rho = (t_{\sf d} {\otimes} t_{\sf d})(\rho) =  (t_{\sf u} {\otimes} t_{\sf u})(\rho)$ it follows that $W({\bf K}_{new}) = W({\bf K})$ and thus $w({\bf K}) = w_{new}({\bf K})$.

Suppose that traversal of the diagram ${\bf K}$ from $P$ to $P_{new}$ in the direction of orientation passes through $m \geq 1$ crossing lines and no local extrema. Observe that $u_{\sf d}(\ell) =  u_{\sf u}(\ell)  = -d$ for all $1 \leq \ell \leq m$; see the discussion preceding Proposition \ref{wTop}. Let $x_1, \ldots, x_{2n}$ be the crossing line decorations. Since $t_{\sf d}^{-d} {\circ} t_{\sf u}^{-d}$ is an algebra automorphism of $A$, and $t_{\sf d}, t_{\sf u}$ commute, we can make the substitution $(t_{\sf d}{\circ} t_{\sf u})^{-d}(x_1\cdots x_m)$ for $(t_{\sf d}^{-d}{\circ} t_{\sf u}^{-d})(x_1) \cdots  (t_{\sf d}^{-d}{\circ} t_{\sf u}^{-d})(x_m).$ Thus 
we have 
$$
W({\bf K}) = (t_{\sf d}{\circ} t_{\sf u})^{-d}(x_1 \cdots  x_m) (t^{u_{\sf d}(m{+}1)}{\circ} t^{u_{\sf u}(m{+}1)})(x_{m+1}) \cdots (t_{\sf d}^{u_{\sf d}(2n)}{\circ} t_{\sf u}^{u_{\sf u}(2n)}) (x_{2n})
$$
and
$$ 
W_{new}({\bf K}) =  (t_{\sf d}^{u_{\sf d}(m{+}1)}{\circ} t_{\sf u}^{u_{\sf u}(m{+}1)})(x_{m+1}) \cdots (t_{\sf d}^{u_{\sf d}(2n)}{\circ} t_{\sf u}^{u_{\sf u}(2n)}) (x_{2n})x_1\cdots x_m 
$$
Since
\begin{eqnarray*}
\lefteqn{{\sf tr}(G^d(t_{\sf d}{\circ} t_{\sf u})^{-d}(a_1 \cdots  a_m) (t_{\sf d}^{u_{\sf d}(m{+}1)}{\circ} t_{\sf u}^{u_{\sf u}(m{+}1)})(a_{m+1}) \cdots (t_{\sf d}^{u_{\sf d}(2n)}{\circ} t_{\sf u}^{u_{\sf u}(2n)}) (a_{2n}))} \\ 
& = & {\sf tr}( G^d(G^{-d}a_1 \cdots  a_mG^d (t_{\sf d}^{u_{\sf d}(m{+}1)}{\circ} t_{\sf u}^{u_{\sf u}(m{+}1)})(a_{m+1}) \cdots (t_{\sf d}^{u_{\sf d}(2n)}{\circ} t_{\sf u}^{u_{\sf u}(2n)}) (a_{2n})) \\
& = & {\sf tr}( a_1 \cdots  a_mG^d (t_{\sf d}^{u_{\sf d}(m{+}1)}{\circ} t_{\sf u}^{u_{\sf d}(m{+}1)})(a_{m+1}) \cdots (t_{\sf d}^{u_{\sf d}(2n)}{\circ} t_{\sf u}^{u_{\sf u}(2n)}) (a_{2n})) \\
& = & {\sf tr}( G^d (t_{\sf d}^{u_{\sf d}(m{+}1)}{\circ} t_{\sf u}^{u_{\sf u}(m{+}1)})(a_{m+1}) \cdots (t_{\sf d}^{u_{\sf d}(2n)}{\circ} t_{\sf u}^{u_{\sf u}(2n)}) (a_{2n})a_1 \cdots  a_m)
\end{eqnarray*}
for all $a_1, \ldots, a_{2n} \in A$ it follows that ${\sf tr}(G^d w({\bf K})) = {\sf tr}(G^d w_{new}({\bf K}))$ in the second case. Thus ${\sf tr}(G^d w({\bf K}))$ does not depend on the starting point in any event and therefore  
\begin{equation}\label{fAcKnots}
{\bf Inv}_{A, \, {\sf tr}}({\bf K}) = {\sf tr}(G^dw({\bf T}))
\end{equation}
is a well-defined for all ${\bf K} \in {\SL Knot}$. We regard (\ref{fAcKnots}) as defining a function from ${\SL Knot}$ to $k$ which we refer to as ${\bf Inv}_{A, \, {\sf tr}}$ by slight abuse of notation.  For the same reasons that $w({\bf T})$ is not affected by regular isotopy moves $w({\bf K})$ unaffected regular isotopy moves since we may assume that the starting point is not in a local area of the diagram ${\bf K}$ under consideration.

Notice that the oriented knot diagram ${\bf K}$ is regularly isotopic to an oriented knot diagram ${\bf K}({\bf T})$, where ${\bf K}({\bf T})$ is one of
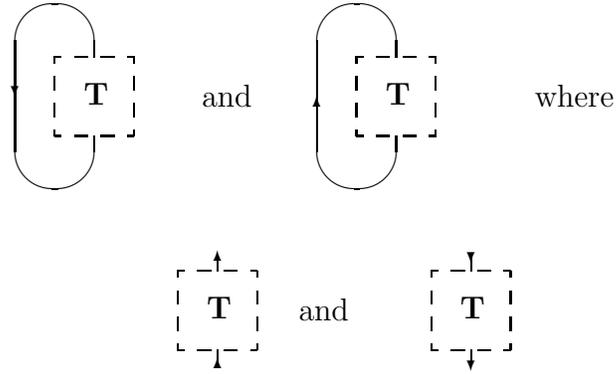
\begin{center}
\mbox{
\begin{picture}(60,80)(-5,-25)
\put(0,-5){\line(0,1){40}}
\put(15,-5){\oval(30, 30)[b]}
\put(15,35){\oval(30, 30)[t]}
\put(30,0){\line(0,-1){5}}
\put(30,30){\line(0,1){5}}
\put(15,0){\dashbox{5}(30,30)}
\put(26,12){$\bf T$}
\put(0,15){\vector(0, -1){0}}
\end{picture}
\raisebox{7ex}{\quad and \quad}
\begin{picture}(60,80)(-5,-25)
\put(0,-5){\line(0,1){40}}
\put(15,-5){\oval(30, 30)[b]}
\put(15,35){\oval(30, 30)[t]}
\put(30,0){\line(0,-1){5}}
\put(30,30){\line(0,1){5}}
\put(15,0){\dashbox{5}(30,30)}
\put(26,12){$\bf T$}
\put(0,15){\vector(0, 1){0}}
\end{picture}}
\raisebox{7ex}{\qquad where \qquad}
\mbox{
\begin{picture}(30,80)(0,-25)
\put(15,0){\line(0,-1){5}}
\put(15,-3){\vector(0, 1){0}}
\put(15,30){\line(0,1){5}}
\put(15,38){\vector(0, 1){0}}
\put(0,0){\dashbox{5}(30,30)}
\put(11,12){$\bf T$}
\end{picture}
\raisebox{7ex}{\quad and \qquad}
\begin{picture}(30,80)(0,-25)
\put(15,0){\line(0,-1){5}}
\put(15,-8){\vector(0, -1){0}}
\put(15,30){\line(0,1){5}}
\put(15,33){\vector(0, -1){0}}
\put(0,0){\dashbox{5}(30,30)}
\put(11,12){$\bf T$}
\end{picture}}
\end{center}
are oriented $1$--$1$ tangle diagrams. Since the Whitney degree is a regular isotopy invariant of oriented knot diagrams, the Whitney degrees of ${\bf K}$ and ${\bf K}({\bf T})$ are the same. 

The reader may have noticed that the definition of ${\bf Inv}_{A, {\sf tr}}({\bf K})$ for ${\bf K} \in {\SL Knot}$ does not require the full use of the axioms for a twist oriented quantum algebra.
\begin{TH}\label{fCcKINV}
Let $(A, \rho, t_{\sf d}, t_{\sf u})$ be an oriented quantum algebra over the field $k$, suppose that $G \in A$ is invertible and satisfies $t_{\sf d} {\circ} t_{\sf u}(a) = GaG^{-1}$ for all $a \in A$,  let ${\sf tr} \in A^*$ be a tracelike element and let ${\bf Inv}_{A, \, {\sf tr}} : {\SL Knot} \la k$ be the function defined by (\ref{fAcKnots}). 
\begin{enumerate}
\item[{\rm a)}]
Suppose that ${\bf K}, {\bf K}' \in {\SL Knot}$ are regularly isotopic. Then ${\bf Inv}_{A, \, {\sf tr}}({\bf K}) = {\bf Inv}_{A, \, {\sf tr}}({\bf K}')$.
\item[{\rm b)}]
Suppose that ${\bf K} \in {\SL Knot}$ and that ${\bf K}$ is regularly isotopic to ${\bf K}({\bf T})$ for some ${\bf T} \in {\SL Tang}$. Then 
$$
{\bf Inv}_{A, \, {\sf tr}}({\bf K}) = {\sf tr}(G^d{\bf w}_A({\bf T})),
$$
where $d$ is the Whitney degree of ${\bf K}$.
\item[{\rm c)}]
${\bf Inv}_{A, \, {\sf tr}}({\bf K^{op}})  = {\bf Inv}_{A^{op}, \, {\sf tr}}({\bf K})$ for all ${\bf K} \in {\SL Knot}$.
\item[{\rm d)}]
Suppose that $A_s$ is the standard oriented quantum algebra associated with $A$. Then ${\bf Inv}_{A_s, \, {\sf tr}}({\bf K}) = {\bf Inv}_{A, \, {\sf tr}}({\bf K})$ for all ${\bf K} \in {\SL Knot}$.
\end{enumerate}
\end{TH}

\pf
We have established parts a) and b). As for part c), we may assume that ${\bf K} = {\bf K}({\bf T})$ for some ${\bf T} \in {\SL Tang}$. Using part b) and Proposition \ref{wTop} we see that  
$$
{\bf Inv}_{A, \, {\sf tr}}({\bf K^{op}}) = {\sf tr}(G^{-d}{\bf w}_A({\bf T^{op}})) = {\sf tr}({\bf w}_{A^{op}}({\bf T})G^{-d}) = {\bf Inv}_{A^{op}, \, {\sf tr}}({\bf K}).
$$
Part d) follows by part b) and Theorem \ref{11Standard}. This concludes our proof.
\qed
\medskip

We now turn to links. Let ${\bf L} \in {\SL Link}$ be an oriented link diagram with components ${\bf L}_1, \ldots, {\bf L}_r$. To construct ${\bf Inv}_{A, \, {\sf tr}}({\bf L})$ we  modify slightly the steps of the construction of ${\bf Inv}_{A, \, {\sf tr}}({\bf K})$, where ${\bf K} \in {\cal K}$.

Decorate the crossings of ${\bf L}$ according to the conventions of Section \ref{secINV11T}. For each $1 \leq \ell \leq r$ let $d_\ell$ denote the Whitney degree of the component ${\bf L}_\ell$ and choose a point $P_\ell$ on a vertical line of ${\bf L}_\ell$. (As in the case of knot diagrams we can always assume that each component of ${\bf L}$ has a vertical line.) Traverse ${\bf L}_\ell$ in the direction of the orientation, beginning at $P_\ell$, and label the crossing lines (if any) contained in ${\bf L}_\ell$ by $(\ell {:}1), (\ell {:}2), \ldots (\ell {:} m)$ in the order encountered. For $1 \leq \imath \leq m$ let $u_{\sf d}(\ell {:} \imath)$ denote the number of local extrema which of type (${\sf d}_+$) minus the number of type (${\sf d}_-$) which are encountered in the part of the traversal of ${\bf L}_\ell$ from the line labeled $\imath$ back to the starting point $P_\ell$. Define $u_{\sf u}(\ell)$ in the same manner, where (${\sf u}_+$) and (${\sf u}_-$) replace (${\sf d}_+$) and (${\sf d}_-$) respectively. Let $x_{(\ell {:} \imath)}$ be the decoration on the line labeled $(\ell {:} \imath)$.

We next define a formal word $W({\bf L}_\ell)$ as follows. If ${\bf L}_\ell$ contains no crossing lines we set $W({\bf L}_\ell) = 1$; otherwise we set 
$$
W({\bf L}_\ell) = t_{\sf d}^{u_{\sf d}(\ell {:} 1)}{\circ} t_{\sf u}^{u_{\sf u}(\ell {:} 1)}(x_{(\ell {:}1)}) \cdots  t_{\sf d}^{u_{\sf d}(\ell {:} m)}{\circ} t_{\sf u}^{u_{\sf u}(\ell {:} m)}(x_{(\ell {:}m)}).
$$
Now let $w({\bf L}_1) {\otimes} \cdots {\otimes}w({\bf L}_r) \in A{\otimes} \cdots {\otimes} A$ be obtained by substituting copies of $\rho$ and $\rho^{-1}$ into the formal tensor $W({\bf L}_{\ell_1}) {\otimes} \cdots {\otimes}W({\bf L}_{\ell_r})$. In light of our discussion ${\bf Inv}_{A, \, {\sf tr}}({\bf K})$, where ${\bf K} \in {\SL Knot}$, it is a small exercise to show that  the scalar
\begin{equation}\label{EqfAtrL}
{\bf Inv}_{A, \, {\sf tr}}({\bf L}) = {\sf tr}(G^{d_1}w({\bf L}_1)) \cdots {\sf tr}(G^{d_r}w({\bf L}_r))
\end{equation}
does not depend on the particular choice of $P_1, \cdots, P_r$ and is not affected by regular isotopy moves. We note that the equation $u_{\sf d}(\imath) = u_{\sf u}(\imath)$ for all $1 \leq \imath \leq m$ holds when the crossing lines in ${\bf L}_\ell$ are directed upward and traversal of ${\bf L}_\ell$ begins in the upward direction; see the discussion preceding Proposition \ref{wTop}.
\begin{TH}\label{fCcLINV}
Suppose that $(A, \rho, t_{\sf t}, t_{\sf u}, G)$ is a twist oriented quantum algebra over the field $k$, let ${\sf tr} \in A^*$ be a $t_{\sf d}^*, t_{\sf u}^*$-invariant tracelike element and let ${\bf Inv}_{A, \, {\sf tr}} : {\SL Link} \la k$ be the function defined by (\ref{EqfAtrL}). 
\begin{enumerate}
\item[{\rm a)}]
Suppose that ${\bf L}, {\bf L}' \in {\SL Link}$ are regularly isotopic. Then ${\bf Inv}_{A, \, {\sf tr}}({\bf L}) = {\bf Inv}_{A, \, {\sf tr}}({\bf L}')$.
\item[{\rm b)}]
Suppose that $A_s$ is the standard twist oriented quantum algebra associated with $A$. Then ${\bf Inv}_{A_s, \, {\sf tr}}({\bf L}) = {\bf Inv}_{A, \, {\sf tr}}({\bf L})$ for all ${\bf L} \in {\SL Link}$.
\end{enumerate} \qed
\end{TH}
\medskip

We next consider the relationship between invariants of the type ${\bf Inv}_{A, \, {\sf tr}}({\bf L})$ and morphisms of twist oriented quantum algebras. We define a {\em morphism} $f : (A, \rho, t_{\sf d}, t_{\sf u}, G) \la (A', R', t'_{\sf d}, t'_{\sf u}, G')$ {\em of twist oriented quantum algebras} to be a morphism $f : (A, \rho, t_{\sf d}, t_{\sf u}) \la (A', R', t'_{\sf d}, t'_{\sf u})$ of oriented quantum algebras which satisfies $f(G) = G'$. 
\begin{PR}
Suppose that $f : (A, \rho, t_{\sf d}, t_{\sf u}, G) \la (A', R', t'_{\sf d}, t'_{\sf u}, G')$ is a morphism of twist oriented quantum algebras over $k$ and that ${\sf tr}' \in A^{'*}$ is a $t'^*_{\sf d}, t'^*_{\sf u}$-invariant tracelike element. Then ${\sf tr} \in A^*$ defined by ${\sf tr} = {\sf tr}'{\circ}f$ is a $t_{\sf d}^*, t_{\sf u}^*$-invariant tracelike element and ${\bf Inv}_{A, \, {\sf tr}}({\bf L}) = {\bf Inv}_{A', \, {\sf tr}'}({\bf L})$ for all ${\bf L} \in {\SL Link}$.
\end{PR}
\qed
\medskip

There are interesting connections between the invariants described in 
\cite{QC}
for quantum algebras and the invariants described here for oriented quantum algebras. Let ${\SL Tang}^u$ (respectively ${\SL Knot}^u$, ${\SL Link}^u$) be the set of unoriented $1$--$1$ tangle (respectively knot, link) diagrams situated with respect to a fixed vertical. When $(A, \rho, s)$ is a quantum algebra then the function ${\sf Inv}_A : {\SL Tang}^u \la A$ defined in 
\cite[Section 6.1]{QC}
determines a regular isotopy invariant of $1$--$1$ tangles since whenever ${\sf T}, {\sf T}' \in {\SL Tang}^u$ are regularly isotopy invariant then ${\sf Inv}_A({\sf T}) = {\sf Inv}_A({\sf T}')$. We note that ${\sf Inv}_A$ was implicitly defined in 
\cite[Section IV]{GAUSS}.

Now suppose that $(A, \rho, s, G)$ is a twist quantum algebra and ${\sf tr} \in A^*$ is an $s^*$-invariant tracelike element. Then the function ${\sf Inv}_{A, \, {\sf tr}} : {\SL Link}^u \la k$ defined in 
\cite[Section 8]{QC}
determines a regular isotopy invariant of knots and links since whenever ${\sf L}, {\sf L}' \in {\SL Link}^u$ are regularly isotopy invariant ${\sf Inv}_{A, \, {\sf tr}}({\sf L}) = {\sf Inv}_{A, \, {\sf tr}}({\sf L}')$. We note that ${\sf Inv}_{A, \, {\sf tr}}$ was explicitely defined in 
\cite[Section IV]{GAUSS}.

The calculations of ${\sf Inv}_A({\sf T})$, where the initial vertical line of ${\sf T}$ is oriented upward, and ${\sf Inv}_{A, \, {\sf tr}}({\sf L})$ are made in the same manner as the calculations of their counterparts ${\bf Inv}_A({\bf T})$ and ${\bf Inv}_{A, \, {\sf tr}}({\bf L})$ for oriented quantum algebras and oriented twist quantum algebras respectively. The decorated crossings 
%
%
\mbox{
\begin{picture}(80,50)(-10,-10)
\put(0,0){\line(1,1){30}}
\put(0,30){\line(1,-1){13}}
\put(30,0){\line(-1,1){13}}
\put(-5, 3){$e \; \bullet$}
\put(22, 3){$\bullet \; e'$}
\put(0, 30){\vector(-1, 1){0}}
\put(30, 30){\vector(1, 1){0}}
\end{picture} }
\quad 
and 
\quad 
%
\mbox{
\begin{picture}(80,50)(-10,-10)
\put(0,0){\line(1,1){13}}
\put(30,0){\line(-1,1){30}}
\put(30,30){\line(-1,-1){13}}
\put(-9, 21){$E \; \bullet$}
\put(22, 21){$\bullet \; E'$}
\put(0, 30){\vector(-1, 1){0}}
\put(30, 30){\vector(1, 1){0}}
\end{picture} }
in oriented diagrams are replaced by 
%
%
\mbox{
\begin{picture}(80,50)(-10,-10)
\put(0,0){\line(1,1){30}}
\put(0,30){\line(1,-1){13}}
\put(30,0){\line(-1,1){13}}
\put(-20, 3){$s(e) \; \bullet$}
\put(22, 3){$\bullet \; e'$}
\end{picture} }
\quad 
and 
\quad 
%
\mbox{
\begin{picture}(80,50)(-10,-10)
\put(0,0){\line(1,1){13}}
\put(30,0){\line(-1,1){30}}
\put(30,30){\line(-1,-1){13}}
\put(-6, 21){$e \; \bullet$}
\put(22, 21){$\bullet \; e'$}
\end{picture} }
respectively. Traversal in the unoriented case begins on a vertical line and proceeds in the upward direction. The rules for traversal of local extrema are the in the unoriented case are the rules in the oriented case where $s$ replaces $t_{\sf d}$ and $t_{\sf u}$. 

Let $u({\cal D})$ denote an oriented diagram ${\cal D}$ with its orientation removed. The relationship between the invariants associated with a quantum algebra $(A, \rho, s)$ and the invariants  associated with the oriented quantum algebra $(A, \rho, 1_A, s^{-2})$ of part a) of Proposition \ref{MinQAoqa} are described in our next result.
\begin{TH}\label{Thmst}
Let $(A, \rho, s)$ be a quantum algebra over $k$. Then: 
\begin{enumerate}
\item[{\rm a)}]
The equations
$$
{\bf Inv}_{(A, \rho, 1_A, s^{-2})}({\bf T}) = {\sf Inv}_{(A, \rho^{-1}, s^{-1})}(u({\bf T}))
$$
and
$$
{\sf Inv}_{(A, \rho, s)}(u({\bf T})) = {\bf Inv}_{(A, \rho^{-1}, 1_A,  s^2)}({\bf T})
$$
hold for all ${\bf T} \in {\SL Tang}$ whose initial vertical line is oriented upward.
\item[{\rm b)}]
Suppose further $(A, \rho, s, G^{-1})$ is a twist quantum algebra and ${\sf tr} \in A^*$ is an $s^*$-invariant tracelike element. Then 
$$
{\bf Inv}_{(A, \rho, 1_A, s^{-2}, G), \, {\sf tr}}({\bf L}) = {\sf Inv}_{(A, \rho^{-1}, s^{-1}, G), \, {\sf tr}}(u({\bf L}))
$$
and
$$
{\sf Inv}_{(A, \rho, s, G^{-1}), \,{\sf tr}}(u({\bf L})) = {\bf Inv}_{(A, \rho^{-1}, 1_A,  s^2, G^{-1}), \,{\sf tr}}({\bf L})
$$ 
for all ${\bf L} \in {\SL Link}$.
\end{enumerate}
\end{TH}

\pf
We need only establish the first equations in parts a) and b). Generally to calculate a regular isotopy invariant of oriented $1$--$1$ tangle, knot or link diagrams we may assume that all crossing lines are directed upward, and we may assume that diagrams have vertical lines oriented in the upward direction. We can assume that traversals begin on such lines.

Write $\rho = e {\otimes}e'$ and $\rho^{-1} = E {\otimes}E'$. Since $\rho^{-1} = (s {\otimes}1_A)(\rho)$ it follows that $e {\otimes}e' = s^{-1}(E) {\otimes} E'$. Thus in the diagrams preceding the statement of the theorem, the oriented crossing decorations representing $e {\otimes}e'$ and $E {\otimes} E'$ are the unoriented crossing decorations representing $s^{-1}(E) {\otimes} E'$ and $E{\otimes}E'$ respectively; that is the decorations associated with the quantum algebra $(A, \rho^{-1}, s^{-1})$.

For a crossing decoration representing $x{\otimes}y$, traversal of the $1$--$1$ tangle or link diagram results in the modification  
$$
t_{\sf d}^{u_{\sf d}} {\circ} t_{\sf u}^{u_{\sf u}}(x){\otimes} t_{\sf d}^{u'_{\sf d}} {\circ} t_{\sf u}^{u'_{\sf u}}(y) = s^{-2u_{\sf u}} {\otimes} s^{-2u'_{\sf u}}
$$
with $(A, \rho, 1_A, s^{-2})$ and results in the modification  
$$
s^{-(u_{\sf d} + u_{\sf u})}(x) {\otimes}s^{-(u'_{\sf d} + u'_{\sf u})}(y) = s^{-2u_{\sf u}}(x) {\otimes} s^{-2u'_{\sf u}}(y)
$$ 
with $(A, \rho^{-1}, s^{-1})$; the last equation holds by our assumptions on the crossings and traversal.
\qed
\medskip

We end this section with a simple example, the oriented Hopf link ${\bf L}_{\rm Hopf}$ depicted below left with components ${\bf L}_1$ and ${\bf L}_2$, reading left to right. The symbol  $\circ$ denotes a starting point for component traversal.
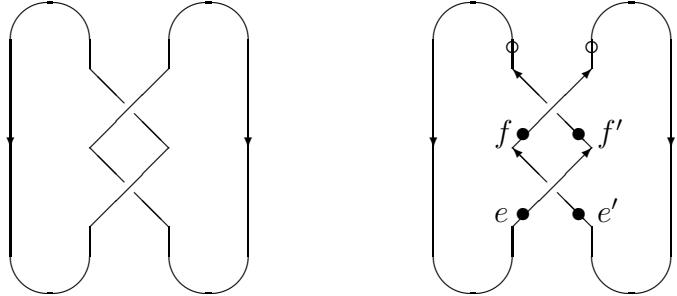
\begin{center}
\begin{picture}(280,140)(0,0)
\put(170,30){\line(0,1){80}}
\put(260,30){\line(0,1){80}}
\put(200,30){\line(0,1){10}}
\put(200,100){\line(0,1){10}}
\put(230,30){\line(0,1){10}}
\put(230,100){\line(0,1){10}}
\put(170,70){\vector(0, -1){0}}
\put(260,70){\vector(0, -1){0}}
\put(185,110){\oval(30,30)[t]}
\put(245,110){\oval(30,30)[t]}
\put(185,30){\oval(30,30)[b]}
\put(245,30){\oval(30,30)[b]}
\put(200,40){\line(1,1){30}}
\put(200,70){\line(1,1){30}}
\put(230,70){\vector(1, 1){0}}
\put(200,70){\line(1,-1){13}}
\put(200,70){\vector(-1, 1){0}}
\put(230,100){\vector(1, 1){0}}
\put(200,100){\vector(-1, 1){0}}
\put(200,100){\line(1,-1){13}}
\put(230,40){\line(-1,1){13}}
\put(230,70){\line(-1,1){13}}
\put(10,30){\line(0,1){80}}
\put(100,30){\line(0,1){80}}
\put(40,30){\line(0,1){10}}
\put(40,100){\line(0,1){10}}
\put(70,30){\line(0,1){10}}
\put(70,100){\line(0,1){10}}
\put(10,70){\vector(0, -1){0}}
\put(100,70){\vector(0, -1){0}}
\put(25,110){\oval(30,30)[t]}
\put(85,110){\oval(30,30)[t]}
\put(25,30){\oval(30,30)[b]}
\put(85,30){\oval(30,30)[b]}
\put(40,40){\line(1,1){30}}
\put(40,70){\line(1,-1){13}}
\put(70,40){\line(-1,1){13}}
\put(40,70){\line(1,1){30}}
\put(40,100){\line(1,-1){13}}
\put(70,70){\line(-1,1){13}}
\put(193, 42){$e$}
\put(201, 42){$\bullet$}
\put(232, 42){$e'$}
\put(222, 42){$\bullet$}
\put(193, 72){$f$}
\put(201, 72){$\bullet$}
\put(232, 72){$f'$}
\put(222, 72){$\bullet$}
\put(197, 105){$\circ$}
\put(227, 105){$\circ$}
\end{picture}
\end{center}

Observe that $d_1 = -1$, $d_2 = 1$ and 
$$
{\bf Inv}_{A, \, {\sf tr}}({\bf L}_{\rm Hopf}) = \sum_{\imath = 1}^r \sum_{\jmath = 1}^r {\sf tr}(G^{-1}a_\imath b_\jmath) {\sf tr}(Gb_\imath a_\jmath).
$$
\section{Oriented Quantum Coalgebras and $T$-Form Structures}\label{secOrientCoalg}
In this very brief section we define the coalgebra counterparts of oriented quantum algebras and twist oriented quantum algebras. We describe a very important connection between the coalgebra structures defined in this section and $T$-form structures
\cite[Section 3]{RKQCS}.
In 
\cite{KRNexT}
a general theory of the coalgebra structures discussed below is developed along the lines of 
\cite{QC}.

Let $(C, \Delta, \epsilon )$ be a coalgebra over the field $k$. We denote the coproduct $\Delta (c) \in C{\otimes} C$ symbolically by $\Delta (c) = c_{(1)} {\otimes} c_{(2)}$. Let $b, b' : C {\times} C \la k$ be bilinear forms. Then $b'$ is an inverse of $b$ if 
$$
b(c_{(1)}, d_{(1)})b'(c_{(2)}, d_{(2)}) = \epsilon (c)\epsilon (d) = 
b'(c_{(1)}, d_{(1)})b(c_{(2)}, d_{(2)})
$$
for all $c, d \in C$. If $b$ has an inverse then its inverse is unique and we denote it by $b^{-1}$. 

A {\em strict oriented quantum coalgebra over} $k$ is a quadruple $(C, b, T_{\sf d}, T_{\sf u})$, where $C$ is a coalgebra over $k$, $b : C {\times} C \la k$ is an invertible bilinear form and $T_{\sf d}, T_{\sf u}$ are commuting coalgebra automorphisms of $C$, such that 
\begin{flushleft}
{\rm (qc.1)} 
$b(c_{(1)}, T_{\sf u}(d_{(2)}))b^{-1}(T_{\sf d}(c_{(2)}), d_{(1)}) =  \epsilon (c)\epsilon (d)$ and 
$\phantom{aaaaaaa} b^{-1}(T_{\sf d}(c_{(1)}), d_{(2)})b(c_{(2)}, T_{\sf u}(d_{(1)})) =  \epsilon (c)\epsilon (d)$,
\vskip1\jot
{\rm (qc.2)} 
$b(c, d) = b(T_{\sf d}(c), T_{\sf d}(d)) = b(T_{\sf u}(c), T_{\sf u}(d))$ and
\vskip1\jot
{\rm (qc.3)}
$b(c_{(1)}, d_{(1)})b(c_{(2)}, e_{(1)})b(d_{(2)}, e_{(2)})
=
b(c_{(2)}, d_{(2)})b(c_{(1)}, e_{(2)})b(d_{(1)}, e_{(1)})$
\end{flushleft}
for all $c, d, e \in C$.

The notion of strict oriented quantum coalgebra is dual to the notion of oriented quantum algebra. To see this, let $A$ be a finite-dimensional algebra over $k$. Then the linear dual $A^*$ is a coalgebra over $k$ where $\epsilon (a^*)$ and $\Delta (a^*)$ are computed for $a^* \in A^*$ as follows:
$$
\epsilon (a^*) = a^*(1) \quad \mbox{and} \quad \Delta (a^*) = \sum_{\imath = 1}^r a_\imath^* {\otimes} b_\imath^* \in A^* {\otimes} A^* = (A {\otimes} A)^*,
$$
where $a^*(ab) = \sum_{\imath = 1}^r a_\imath^*(a)b_\imath^*(b)$ for all $a, b \in A^*$. 

Let $\rho \in A{\otimes}A$, suppose that $b : A {\times} A \la k$ is the bilinear form defined by $b(a^*, b^*) = (a^*{\otimes}b^*)(\rho)$ for all $a^*, b^* \in A^*$ and let $t_{\sf d}, t_{\sf u}$ be linear automorphisms of $A$. It is a straightforward exercise to show that $(A^*, b, t_{\sf d}^*, t_{\sf u}^*)$ is a strict oriented quantum coalgebra over $k$ if and only if $(A, \rho, t_{\sf d}, t_{\sf u})$ is an oriented quantum algebra over $k$.

An {\em oriented quantum coalgebra over} $k$ is a quadruple $(C, b, T_{\sf d}, T_{\sf u})$, where $C$ is a coalgebra over $k$, $b : C {\times} C \la k$ is an invertible bilinear form and $T_{\sf d}, T_{\sf u}$ commuting coalgebra automorphisms of $C$ with respect to $\{ b, b^{-1}\}$, such that (qc.1)--(qc.3) hold. A linear automorphism $T$ of $C$ is a coalgebra automorphism with respect to $\{ b, b^{-1}\}$ if $\epsilon {\circ} T = \epsilon$, 
$$
b'(T(c_{(1)}), d)b''(T(c_{(2)}), e) = b'(T(c)_{(1)}, d)b''(T(c)_{(2)}, e)  \;\; \mbox{and}
$$
$$
b'(d, T(c_{(1)}))b''(e, T(c_{(2)})) = b'(d, T(c)_{(1)})b''(e, T(c)_{(2)}) 
$$
for all $b', b'' \in \{ b, b^{-1}\}$ and $c, d, e \in C$. An oriented quantum coalgebra $(C, b, T_{\sf d}, T_{\sf u})$ over $k$ is {\em standard} if $T_{\sf d} = 1_C$ and is {\em balanced} if $T_{\sf d} = T_{\sf u}$, in which case we write $(C, b, T)$ for $(C, b, T_{\sf d}, T_{\sf u})$, where $T = T_{\sf d} = T_{\sf u}$.

Suppose that $(C, b, T)$ is a balanced oriented quantum coalgebra over $k$. Since $T$ is a coalgebra automorphism with respect to $\{ b, b^{-1}\}$ and $b(c, d) = b(T(c), T(d))$ for all $c, d \in C$, it follows that  $b^{-1}(c, d) = b^{-1}(T(c), T(d))$ for all $c, d \in C$ and that
$T^{-1}$ is a coalgebra automorphism with respect to $\{ b, b^{-1}\}$. Thus (qc.1) may be reformulated 
$$
b(T^{-2}(c_{(1)}), d_{(2)})b^{-1}(c_{(2)}, d_{(1)}) = \epsilon (c)\epsilon (d) = 
b^{-1}(c_{(1)}, T^{-2}(d_{(2)}))b(c_{(2)}, d_{(1)}) 
$$
for all $c, d \in C$; that is $(C, b, T^{-2})$ is a $T^{-2}$-form structure. More precisely:
\begin{PR}\label{TMinus2Form}
Let $C$ be a coalgebra over the field $k$, let $b :  C{\times} C \la k$ be an invertible bilinear form, and suppose that $T$ is a coalgebra automorphism of $C$ with respect to $\{ b, b^{-1}\}$. Then the following are equivalent:
\begin{enumerate}
\item[{\rm a)}]
$(C, b, T)$ is a balanced oriented quantum coalgebra.
\item[{\rm b)}]
$(C, b, T^{-2})$ is $T^{-2}$-form structure over $k$. 
\end{enumerate} 
\qed
\end{PR}
\medskip

This proposition establishes a fundamental connection between regular isotopy invariants of oriented knots and links computed by the methods of Section \ref{secInv} and regular isotopy invariants of unoriented knots and links computed by the methods of 
\cite{RKQCS,RADKAUFFinv,QC}.
The connection is discussed in great detail in 
\cite{KRNexT}.

We end this section with the definition of twist oriented quantum coalgebra, a structure which is the counterpart of  twist oriented quantum algebra. For a coalgebra $C$ over $k$ recall that the linear dual $C^*$ is an algebra over $k$.

A {\em twist oriented quantum coalgebra over} $k$ is a quintuple $(C, b, T_{\sf d}, T_{\sf u}, G)$, where $(C, b, T_{\sf d}, T_{\sf u})$ is a strict oriented quantum coalgebra over $k$, $G \in C^*$ is invertible and 
$$
T_{\sf d}^*(G) = T_{\sf u}^*(G) = G  \quad \mbox{and} \quad  T_{\sf d}{\circ} T_{\sf u}(x) = G^{-1}{\rh}x{\lh}G = (G^{-1}{\rh}x){\lh}G
$$ 
for all $x \in C$, where $c^*{\rh}c = c_{(1)}{<}c^*, c_{(2)}{>}$ and $c{\lh}c^* =  {<}c^*, c_{(1)}{>}c_{(2)}$ for all $c^* \in C^*$ and $c \in C$. When the underlying oriented quantum coalgebras structure of a twist oriented quantum coalgebra $(C, b, T_{\sf d}, T_{\sf u} , G)$ over $k$ is balanced we shall write $(C, b, T, G)$ for $(C, b, T_{\sf d}, T_{\sf u}, G)$, where $T = T_{\sf d} = T_{\sf u}$, and call $(C, b, T, G)$ a {\em twist balanced oriented quantum coalgebra over} $k$.

Twist oriented quantum coalgebras over $k$ give rise to regular isotopy invariants of oriented knots and links
\cite{KRNexT}.
As one might suspect, the notion of twist oriented quantum coalgebra is dual to the notion of twist oriented quantum algebra. 
\section{Parameterized Families of Oriented Quantum Algebra Structures on ${\rm M}_n(k)$}\label{secPF}
In this section we study the balanced oriented quantum algebra structures $(A, \rho, t)$ on $A = {\rm M}_n(k)$ where $t$ is the automorphism of Example \ref{ExamMnkRaBC}. The description of $\rho$  involves an ordering on $\{ 1, \ldots, n\}$ which is usually not the standard one. Thus for conceptual reasons at the outset we will regard $A = {\rm M}_S(k)$ as the $k$-algebra with basis of symbols $\{ E_{\imath \, \jmath} \}_{\imath, \jmath \in S}$ which satisfy $E_{\imath \, \jmath}E_{\ell \, m} = \delta_{\jmath \, \ell} E_{\imath \, m}$ for all $\imath, \jmath, \ell, m \in S$, where $S$ is an $n$-element set with no a priori ordering.

Let $\omega_\imath \in k^\star$ for all $\imath \in S$. The linear automorphism $t$ of $A$ determined by 
$$
t(E_{\imath \, \jmath}) = \left(\frac{\omega_\imath}{\omega_\jmath}\right)E_{\imath \, \jmath} 
$$
for all $\imath, \jmath \in S$ is an algebra automorphism. Let $R \in A{\otimes} A$ and write 
$$
\rho = \sum_{\imath, \jmath, \ell, m \in S} \rho_{\imath \, \ell  \, \jmath \, m} E_{\imath \, \jmath} {\otimes} E_{\ell \, m},
$$
where $\rho_{\imath \, \ell  \, \jmath \, m} \in k$. Observe that $\rho = (t {\otimes}t)(\rho)$ if and only if 
$$
\rho_{\imath \, \ell  \, \jmath \, m} = \left(\frac{\omega_\imath \omega_\ell}{\omega_\jmath \, \omega_m}\right) \rho_{\imath \, \ell  \, \jmath \,m}
$$ 
for all $\imath, \jmath, \ell, m \in S$. Thus (qa.2) holds for $\rho$ and $t$ if  
\begin{equation}\label{EqRijlmNot0}
\rho_{\imath \, \ell  \, \jmath \, m} \neq 0 \quad \mbox{implies} \quad \{ \imath, \ell \} = \{ \jmath, m\}.
\end{equation}
We note that (qa.2) is equivalent to (\ref{EqRijlmNot0}) when the $\omega_\imath$'s are algebraically independent over the prime field of $k$. Notice that (\ref{EqRijlmNot0}) is satisfied in Example \ref{ExamMnkRaBC}.

Suppose that (\ref{EqRijlmNot0}) holds. Then $(A, \rho, t)$ is a balanced oriented quantum algebra if and only if 
\begin{equation}\label{EqRilil}
\rho_{\imath \, \ell \, \imath \, \ell} \neq 0 \quad \mbox{for all $\imath, \ell \in S$,} 
\end{equation}
\begin{equation}\label{EqRilli}
\rho_{\imath \, \ell  \, \ell \, \imath} = 0  \quad \mbox{or} \quad \rho_{\ell \, \imath \, \imath \, \ell} = 0  \quad \mbox{for all distinct $\imath, \ell \in S$,} 
\end{equation}
\begin{equation}\label{EqRiuui}
\frac{\rho_{\imath \, u \, u \, \imath}}{\rho_{\imath \, \imath  \, \imath \, \imath}} = \delta_{\imath \, u} + \sum_{\ell \neq \imath} \left(\frac{\rho_{\imath \, \ell \, \ell \, \imath}\rho_{\ell \, u \, u \, \ell}} {\rho_{\imath \, \ell \, \imath \, \ell}\rho_{\ell \, \imath \, \ell \, \imath}}\right)\left(\frac{\omega_\imath}{\omega_\ell }\right)^2 \quad \mbox{for all $\imath, u \in S$,} 
\end{equation}
\begin{equation}\label{EqQYBeS}
\sum_{u, v, w \in T} \rho_{\imath \, k \, u \, v} \rho_{u \, p \, \jmath \, w} \rho_{v \, w \, \ell \, q} 
 = \sum_{u, v, w \in T} \rho_{u \, v \, \jmath \, \ell} \rho_{\imath \, w \, u \, q} \rho_{k \, p \, v \, w} 
\end{equation}
for all $\imath, k, p, \jmath, \ell, q \in T$, where $T$ is any $2$-element subset of $S$, and for all distinct $\imath, \ell, k \in S$:
\begin{equation}\label{EqQYBeS1}
\rho_{\imath \, k \, k \, \imath} \rho_{\imath \, \ell \, \ell \, \imath} = 
\rho_{\imath \, \ell \, \ell \, \imath} \rho_{\ell \, k \, k \, \ell} +
\rho_{\imath \, k \, k \, \imath} \rho_{k \, \ell \, \ell \, k},
\end{equation}  
\begin{equation}\label{EqQYBeS2}
\rho_{\imath \, k \, k \, \imath} \rho_{\ell \, k \, k \, \ell} = 
\rho_{\imath \, \ell \, \ell \, \imath} \rho_{\ell \, k \, k \, \ell} +
\rho_{\imath \, k \, k \, \imath} \rho_{\ell \, \imath \, \imath \, \ell}
\end{equation}  
and 
\begin{equation}\label{EqQYBeS3}
(\rho_{\imath \, k \, k \, \imath})^2 \rho_{k \, \ell \, \ell \, k} + 
\rho_{\imath \, k \, \imath \, k}\rho_{k \, \imath \, k \, \imath} \rho_{\imath \, \ell \, \ell \, \imath} =
(\rho_{k \, \ell \, \ell \, k})^2 \rho_{\imath \, k \, k \, \imath} + 
\rho_{k \, \ell \, k \, \ell}\rho_{\ell \, k \, \ell \, k} \rho_{\imath \, \ell \, \ell \, \imath}.
\end{equation}  

A rather tedious calculation shows that the set of equations (\ref{EqQYBeS})--(\ref{EqQYBeS3}) is equivalent to (qa.3) under the assumption that (\ref{EqRijlmNot0}) and (\ref{EqRilil}) hold. By passing to the dual coalgebra ${\rm C}_n(k) = {\rm M}_n(k)^*$, one sees that this equivalence is established in  
\cite[Lemma 4]{RBracket}. 
Using Proposition \ref{TMinus2Form} and retracing the proof of
\cite[Theorem 2]{RBracket}
we can conclude that the set of statements (\ref{EqRilil})--(\ref{EqRiuui}) is equivalent to $\rho$ is invertible and (qa.1) holds for $\rho$ and $t$ under the assumption that (\ref{EqRijlmNot0}) holds. For the sake of completeness we will sketch a proof that  (\ref{EqRilil})--(\ref{EqRiuui}) are collectively equivalent to $\rho$ is invertible, (qa.1) holds for $\rho$ and $t$ under the assumption that (\ref{EqRijlmNot0}) holds.

Suppose that $\rho$ is invertible and (\ref{EqRijlmNot0}) holds. Let $Q = \rho^{-1}$ and write $Q = \sum_{\imath, \jmath, \ell, m \in S} Q_{\imath \, \ell \, \jmath \, m} E_{\imath \, \jmath} {\otimes} E_{\ell \, m}$ where $Q_{\imath \, \ell \, \jmath \, m}  \in k$. Since $A {\otimes} A$ is a finite-dimensional algebra over $k$ and the set of elements of $A {\otimes} A$ which satisfy (\ref{EqRijlmNot0}) is a subalgebra of $A {\otimes} A$, it follows that $Q_{\imath \, \ell \, \jmath \, m} \neq 0$ implies $\{ \imath, \ell \} = \{ \jmath, m\}$.  Also since $A{\otimes} A$ is finite-dimensional, (qa.1) for $t$ and $\rho$ is equivalent to 
$$
((t{\otimes}1_A)(\rho^{-1}))((1_A {\otimes}t)(\rho)) = 1 {\otimes}1
$$ 
in $A {\otimes} A^{op}$. This equation can be expressed as 
\begin{equation}\label{Eqqa1IN}
\sum_{\jmath, \ell \in S} \left(\frac{\omega_\imath}{\omega_\jmath}\right)\left( \frac{\omega_v}{\omega_\ell}\right) Q_{\imath \, \ell \, \jmath \, m} \rho_{\jmath \, v \, u \, \ell} = \delta_{\imath \, u} \delta_{v \, m}
\end{equation}
for all $\imath, u, v, m \in S$.

Suppose that $\rho$ is invertible, (\ref{EqRijlmNot0}) and (\ref{Eqqa1IN}) hold. We will show that (\ref{EqRilil})--(\ref{EqRiuui}) follow by considering the cases $\imath \neq m$ and $\imath = m$.

Assume first of all that $\imath \neq m$. Since $\{\imath, \ell \} = \{ \jmath, m\}$ if and only if $\imath = \jmath$ and $\ell = m$, in this case (\ref{Eqqa1IN}) boils down to $Q_{\imath \, m \, \imath \, m}\rho_{\imath \, v \, u \, m} = \delta_{\imath \, u} \delta_{v \, m}$ for all $u, v \in S$. This equation holds if and only if $Q_{\imath \, m \, \imath \, m}\rho_{\imath \, m \, \imath \, m} = 1$.

Assume that $\imath = m$. Then (\ref{Eqqa1IN}) is
$$
\sum_{\ell \in S}\left(\frac{\omega_\imath}{\omega_\ell}\right)^2 Q_{\imath \, \ell \, \ell \, \imath} \rho_{\ell \, v \, u \, \ell} =  \delta_{\imath \, u} \delta_{v \, \imath}
$$
for all $\imath, u, v \in S$. If $u \neq v$ then both sides of this equation are zero. Thus (\ref{Eqqa1IN}) is equivalent to 
\begin{equation}\label{Eqqa1IN2}
\sum_{\ell \in S}\left(\frac{\omega_\imath}{\omega_\ell}\right)^2 Q_{\imath \, \ell \, \ell \, \imath} \rho_{\ell \, u \, u \, \ell} = \delta_{\imath \, u}
\end{equation}
for all $u \in S$ when $\imath = m$.

To complete our analysis of the case $\imath = m$ we examine what it means for $Q$ and $\rho$ to be inverses in light of the relation $Q_{\imath \, \ell \, \imath \, \ell} \rho_{\imath \, \ell \, \imath \, \ell} = 1$ for distinct $\imath, \ell \in S$. Since (\ref{EqRijlmNot0}) holds, $Q$ and $\rho$ are inverses if and only if  $Q_{\imath \, \imath \, \imath \, \imath} \rho_{\imath \, \imath \, \imath \, \imath} = 1$ for all $\imath \in S$ and the matrices 
$$
\left(\begin{array}{cc}
\rho_{\imath \, \ell \, \imath \, \ell} & \rho_{\imath \, \ell \, \ell \, \imath} \\
\rho_{\ell \, \imath \, \imath \, \ell} & \rho_{\ell \, \imath \, \ell \, \imath} 
\end{array}\right) 
\quad \mbox{and} \quad 
\left(\begin{array}{cc}
Q_{\imath \, \ell \, \imath \, \ell} & Q_{\imath \, \ell \, \ell \, \imath} \\
Q_{\ell \, \imath \, \imath \, \ell} & Q_{\ell \, \imath \, \ell \, \imath} 
\end{array}\right) 
$$
are inverses when $\imath, \ell \in S$ are distinct.

Suppose that $\imath \neq \ell$ and let $d$ be the determinant of the matrix on the left above. Then $d \neq 0$, and the calculation
$$
\rho_{\imath \, \ell \, \imath \, \ell} \rho_{\ell \, \imath \, \ell \, \imath} - \rho_{\imath \, \ell \, \ell \, \imath}\rho_{\ell \, \imath \, \imath \, \ell} = d = dQ_{\imath \, \ell \, \imath \, \ell} \rho_{\imath \, \ell \, \imath \, \ell} = \rho_{\ell \, \imath \, \ell \, \imath} \rho_{\imath \, \ell \, \imath \, \ell}
$$
shows that $\rho_{\imath \, \ell \, \ell \, \imath}\rho_{\ell \, \imath \, \imath \, \ell} = 0$. Since $dQ_{\imath \, \ell \, \ell \, \imath} = - \rho_{\imath \, \ell \, \ell \, \imath}$ we have 
$$
Q_{\imath \, \ell \, \ell \, \imath} = - \rho_{\imath \, \ell \, \ell \, \imath}/(\rho_{\imath \, \ell \, \imath \, \ell} \rho_{\ell \, \imath \, \ell \, \imath}).
$$
At this point (\ref{EqRiuui}) is easily deduced from (\ref{Eqqa1IN2}).

We have shown that $\rho$ invertible, (\ref{EqRijlmNot0}) and (qa.1) for $t$ and $\rho$ imply (\ref{EqRilil})--(\ref{EqRiuui}). It is a straightforward exercise to show that (\ref{EqRijlmNot0}) and (\ref{EqRilil}) imply $\rho$ is invertible. Now it is easy to see that (\ref{EqRijlmNot0})--(\ref{EqRiuui}) imply that $\rho$ is invertible and that (qa.1) holds for $t$ and $\rho$.

Before continuing we record a description of $\rho^{-1} = Q$:
$$
Q_{\imath \, \ell \, \jmath \, m} \neq 0 \qquad \mbox{implies} \qquad \{ \imath, \ell \} = \{ \jmath, m\},
$$
$$
Q_{\imath \, \ell \, \imath \, \ell} = \frac{1}{\rho_{\imath \, \ell \, \imath \, \ell}} \qquad \mbox{
for all $\imath, \ell \in S$,}
$$ 
and 
$$
Q_{\imath \, \ell \, \ell \,\imath} = -\frac{\rho_{\imath \, \ell \, \ell \, \imath}}{\rho_{\imath \, \ell \, \imath \, \ell}\rho_{\ell \, \imath \, \ell \, \imath}} \qquad \mbox{when $\imath, \ell \in S$ are distinct.}
$$ 

We now determine all invertible $R \in A {\otimes} A$ which satisfy (\ref{EqRiuui})--(\ref{EqQYBeS3}), given that (\ref{EqRijlmNot0})--(\ref{EqRilli}) hold. 

Suppose that (\ref{EqRijlmNot0})--(\ref{EqRilli}) hold. We have noted that (\ref{EqRijlmNot0}) and (\ref{EqRilil}) imply $\rho$ is invertible. Now $\rho_{\imath \, \ell \, \imath \, \ell} \neq 0$ for all $\imath, \ell \in S$, and all of the other $\rho_{\imath \, \ell \, \jmath \, m}$'s are zero with the possible exception of coefficients of the form $\rho_{\imath \, \ell \, \ell \, \imath}$ where $\imath \neq \ell$. Our analysis of $\rho$ is based on whether or not $\rho_{\imath \, \ell \, \ell \, \imath}$ is zero.

We define an ordering on $S$ as follows:
$$
\imath \prec \ell \qquad \mbox{if and only if} \qquad  \imath \neq \ell \;\; \mbox{and} \;\; \rho_{\imath \, \ell \, \ell \, \imath} \neq 0.
$$
Observe that the pair of statements (\ref{EqQYBeS1}) and (\ref{EqQYBeS2}) is equivalent to the pair of statements:
\begin{equation}\label{EqOrder1}
\imath \prec \ell \;\; \mbox{and} \;\; \ell \prec k \qquad \mbox{implies} \qquad \imath \prec k 
\;\; \mbox{and} \;\; \rho_{\imath \, \ell \, \ell \, \imath} = \rho_{\ell \, k \, k \, \ell} = \rho_{\imath \, k \, k \, \imath}
\end{equation}
and 
\begin{equation}\label{EqOrder2}
\imath \prec \ell \;\; \mbox{and} \;\; \imath \prec k, \;\; \mbox{or} \;\; 
\ell \prec \imath \;\; \mbox{and} \;\; k \prec \imath, 
\quad \mbox{implies} \quad \ell \prec k \;\; \mbox{or} \;\; k \prec \ell
\end{equation}
when $\imath, \ell, k \in S$ are distinct. Given (\ref{EqOrder1}), and hence given (\ref{EqOrder1}) and (\ref{EqOrder2}), observe that (\ref{EqQYBeS3}) is equivalent to 
\begin{equation}\label{EqOrder3}
\imath \prec \ell \qquad \mbox{implies} \qquad \rho_{\imath \, k \, \imath \, k} \rho_{k \, \imath \, k \, \imath} = \rho_{\ell \, k \, \ell \, k} \rho_{k \, \ell \, k \, \ell}
\end{equation}
when $\imath, \ell, k \in S$ are distinct.

We will call a subset $T$ of $S$ an $\rho$-{\em component} if (a) for any two $\imath, \jmath \in T$ either $\imath \prec \jmath$ or $\jmath \prec \imath$ and (b) if $T'$ is a subset of $S$ which satisfies (a) and $T \subseteq T'$ then $T = T'$. It is clear that every element of $S$ is contained in an $\rho$-component. 

Suppose that (\ref{EqQYBeS1})--(\ref{EqQYBeS3}) hold, or equivalently that (\ref{EqOrder1})--(\ref{EqOrder3}) hold, in addition to (\ref{EqRijlmNot0})--(\ref{EqRilli}). By (\ref{EqOrder1}) and (\ref{EqOrder2}) if $T, T'$ are $\rho$-components then either $T = T'$ or $T{\cap}T' = \emptyset$. Thus the $\rho$-components partition $S$.

Observe that each $\rho$-component is well-ordered by $\imath \preceq \jmath$ if and only if $\imath = \jmath$ or $\imath \prec \jmath$. Note that (\ref{EqOrder1}) and (\ref{EqOrder2}) also imply that if $\imath$ and $\jmath$ belong to different $\rho$-components then $\imath \not \prec \jmath$ and $\jmath \not \prec \imath$. Thus $\rho_{\imath \, \jmath \, \jmath \, \imath} = 0$ if $\imath$ and $\jmath$ belong to different $\rho$-components, and exactly one of $\rho_{\imath \, \jmath \, \jmath \, \imath}$ and $\rho_{\jmath \, \imath \, \imath \, \jmath}$ is zero for distinct $\imath, \jmath$ which belong to the same $\rho$-component.

We now examine what it means for (\ref{EqQYBeS}) to hold when (\ref{EqRijlmNot0})--(\ref{EqRilli}) and (\ref{EqQYBeS1})--(\ref{EqQYBeS3}) hold. Let $T$ be a two-element subset of $S$. If the elements of $T$ belong to different $\rho$-components then (\ref{EqQYBeS}) holds. Consequently we need only consider (\ref{EqQYBeS}) when $T$ lies in an $\rho$-component which therefore has at least two elements. 

Let ${\cal S}$ be an $\rho$-component with at least two elements. By virtue of (\ref{EqOrder3}) the $\rho_{\imath \, \ell \, \imath \, \ell}\rho_{\ell \, \imath \, \ell \, \imath}$'s for all $\imath, \ell \in {\cal S}$ distinct have the same value, which we denote by $b\!c$.

Suppose that $T = \{ \imath, \jmath \} \subseteq {\cal S}$ and has two elements. We may assume that $\imath \prec \jmath$. Then (\ref{EqQYBeS}) holds for $T$  if and only if 
$$
\left(\begin{array}{cccc}
\rho_{\imath \, \imath \, \imath \, \imath} & 0 & 0 & 0 \\
0 & \rho_{\imath \, \jmath \, \imath \, \jmath} & \rho_{\imath \, \jmath \, \jmath \, \imath} & 0 \\
0 & 0 & \rho_{\jmath \, \imath \, \jmath \, \imath} & 0 \\
0 & 0 & 0 & \rho_{\jmath \, \jmath \, \jmath \, \jmath} 
\end{array}\right)
$$
satisfies the quantum Yang--Baxter equation. This matrix has the form
$$
\left(\begin{array}{cccc}
a & 0 & 0 & 0 \\
0 & b & x & 0 \\
0 & 0 &c & 0 \\
0 & 0 & 0 & d 
\end{array}\right)
$$
where $a, b, c, d, x \in k^\star$. A straightforward calculation shows that the preceding matrix satisfies the quantum Yang--Baxter equation if and only if  $x(a^2 - bc) = x^2a$ and $x(d^2 - bc) = x^2d$, or equivalently
$$
a^2 \neq bc, \quad a = d \;\; \mbox{or} \;\; ad = -bc, \;\; x = a - bc/a,
$$
in which case  $x = d - bc/d$ also. At this point we simplify notation by writing $a_\ell = \rho_{\ell \, \ell \, \ell \, \ell}$ for $\ell \in {\cal S}$. Observe that $\rho_{\imath \, \jmath \, \jmath \, \imath} = a_\imath - b\!c/a_\imath  = a_\jmath - b\!c/a_\jmath$ when (\ref{EqQYBeS}) holds. In this case all such $\rho_{\imath \, \jmath \, \jmath \, \imath}$'s have the same value. We have shown that (\ref{EqQYBeS}) holds for all two element subsets $T \subseteq {\cal S}$ if and only if there are $b\!c, x \in k^\star$ such that $a_\imath^2 \neq b\!c$ for all $\imath \in {\cal S}$, and 
\begin{equation}\label{Eq8Two}
a_\imath = a_\jmath \;\; \mbox{or} \;\; a_\imath a_\jmath = -b\!c, \quad
\rho_{\imath \, \jmath \, \jmath \, \imath} = a_\imath - b\!c/a_\imath  = a_\jmath - b\!c/a_\jmath = x
\end{equation}
for all $\imath, \jmath  \in {\cal S}$ such that $\imath \prec \jmath$.

Assume that (\ref{EqRijlmNot0})--(\ref{EqRilli}) and (\ref{EqQYBeS})--(\ref{EqQYBeS3}) hold, or equivalently (\ref{EqRijlmNot0})--(\ref{EqRilli}), (\ref{EqQYBeS}) and (\ref{EqOrder1})--(\ref{EqOrder3}) hold. We will find a necessary and sufficient condition for (\ref{EqRiuui}) to hold for all $\imath, u \in {\cal S}$.

Observe that (\ref{EqRiuui}) holds when $u = \imath$. If $u \neq \imath$ and $\imath \not \prec u$ or $u \prec \imath$ then both sides of (\ref{EqRiuui}) is zero by (\ref{EqOrder1}). Consequently we need only consider (\ref{EqRiuui}) when $\imath \prec u$; in particular when $\imath, u$ in the same $\rho$-component which we may assume is ${\cal S}$.

Suppose that $\imath, u \in {\cal S}$ and $\imath \prec u$. Then (\ref{EqRiuui}) is equivalent to 
\begin{equation}\label{EqA}
\frac{1}{x}\left( a_\imath - a_u\left(\frac{\omega_\imath}{\omega_u}\right)^2\right) = 1 + \sum_{\imath \prec \ell \prec u}\left(\frac{\omega_\imath}{\omega_\ell}\right)^2.
\end{equation}
If $\imath$ is the immediate predecessor of $u$ in the well-ordering $\preceq$ on ${\cal S}$ then (\ref{EqA}) is 
$$
\frac{1}{x}\left( a_\imath - a_u\left(\frac{\omega_\imath}{\omega_u}\right)^2\right) = 1,
$$
or
$$
a_\imath - a_u\left(\frac{\omega_\imath}{\omega_u}\right)^2 = x = a_\imath - \frac{b\!c}{a_\imath}.
$$
Thus (\ref{EqA}) is equivalent to $(\omega_\imath/\omega_u)^2 = b\!c/(a_\imath a_u)$ when $\imath$ is the immediate predecessor of $u$. Generally, write $\imath = \imath_0 \prec \imath_1 \prec \cdots \prec \imath_m = u$ where $\imath_{r-1}$ is the immediate predecessor of $\imath_r$ for $1 \leq r \leq m$. Since $(\omega_\imath/\omega_u)^2 = (\omega_{\imath_0}/\omega_{\imath_1})^2 (\omega_{\imath_1}/\omega_{\imath_2})^2 \cdots (\omega_{\imath_{r-1}}/\omega_{\imath_r})^2$ we have
\begin{equation}\label{EqB}
\left(\frac{\omega_\imath}{\omega_u}\right)^2 = \frac{b\!c}{a_\imath a_u}\left( \prod_{\imath \prec \jmath \prec u}\frac{b\!c}{a^2_\jmath}\right).
\end{equation}
Let $e$ be smallest element of ${\cal S}$ with respect to the well-ordering $\preceq$. Then (\ref{EqB}) is equivalent to 
\begin{equation}\label{EqC}
\omega_u^2 =  \frac{a_ea_u}{b\!c} \left(\prod_{e \prec \jmath \prec u} \frac{a_\jmath^2}{b\!c}\right) \omega_e^2
\end{equation}
for all $u \in {\cal S}{\setminus}e$.

We have shown that (\ref{EqRiuui}) implies (\ref{EqC}) under the assumption that (\ref{EqRijlmNot0})--(\ref{EqRilli}) and (\ref{EqQYBeS})--(\ref{EqQYBeS3}) hold. Under this assumption we will show that (\ref{EqC}) implies (\ref{EqRiuui}).

Suppose that (\ref{EqC}) holds for $u, \imath \in {\cal S}$. As we have noted, to show that (\ref{EqRiuui}) holds we may assume that $\imath, u \in {\cal S}$ and that $\imath \prec u$. Let $x_\jmath = b\!c/a_\jmath^2$ for $\jmath \in {\cal S}$. Then using (\ref{EqB}) we see that (\ref{EqA}), the equivalent to (\ref{EqRiuui}) in this case, can be formulated as
$$
\frac{1}{x}\left(a_\imath - a_u\left(\frac{b\!c}{a_\imath a_u} \prod_{\imath \prec \jmath \prec u} x_\jmath \right)\right) = 1 + \sum_{\imath \prec \ell \prec u} \frac{b\!c}{a_\imath a_\ell} \left( \prod_{\imath \prec \jmath \prec \ell}x_\jmath \right),
$$
and thus can be written
$$
 \frac{a_\imath}{x}\left(1 -  \prod_{\imath \preceq \jmath \prec u} x_\jmath \right) = 1 + \sum_{\imath \prec \ell \prec u} \frac{b\!c}{a_\imath a_\ell}\left( \prod_{\imath \prec \jmath \prec \ell}x_\jmath \right),
$$
which in turn, using (\ref{Eq8Two}), can be formulated as
$$
1 -  \prod_{\imath \preceq \jmath \prec u} x_\jmath = 
\left(a_\imath - \frac{b\!c}{a_\imath}\right)\frac{1}{a_\imath} + \sum_{\imath \prec \ell \prec u} \left(\frac{b\!c}{a_\imath a_\ell} \left( \prod_{\imath \prec \jmath \prec \ell}x_\jmath \right) \right)\left(a_\ell - \frac{b\!c}{a_\ell}\right)\frac{1}{a_\imath},
$$
which can be written 
$$
\left(\prod_{\imath \preceq \jmath \prec u} x_\jmath\right) - x_\imath = \sum_{\imath \prec \ell \prec u}\left( \prod_{\imath \preceq \jmath \prec \ell} x_\jmath \right)\left(x_\ell - 1\right).
$$
The last equation is always true. For suppose $z_1, \ldots, z_m \in k$ and $1 \leq \imath < m$. Then
$$
\left(\prod_{\imath \leq \jmath < m} z_\jmath\right) - z_\imath = \sum_{\imath < \ell < m}\left( \prod_{\imath \leq \jmath < \ell} z_\jmath \right)\left(z_\ell - 1\right)
$$
follows by induction on $m - \imath$. We have proved:
\begin{TH}\label{ThmLast}
Let $A = {\rm M}_n(k)$ be the algebra of $n{\times}n$ matrices over the field $k$. Let $\omega_1, \ldots, \omega_n \in k^\star$ and $t$ be the algebra automorphism of $A$ determined by 
$$
t(E_{\imath \, \jmath} ) = \left(\frac{\omega_\imath}{\omega_\jmath}\right ) E_{\imath \, \jmath}
$$ 
for all $1 \leq \imath, \jmath \leq n$, where $\{ E_{\imath \, \jmath} \}_{1 \leq \imath, \jmath  \leq n}$ is the standard basis for $A$. Let $R \in A{\otimes} A$ and write 
$$
\rho = \sum_{\imath, \jmath, \ell, m = 1}^n \rho_{\imath \, \ell \, \jmath \, m} E_{\imath \, \jmath} {\otimes} E_{\ell \, m}
$$ 
where $\rho_{\imath \, \ell \, \jmath \, m} \in k$. Assume that $\rho_{\imath \, \ell \, \jmath \, m}  \neq 0$ implies $\{ \imath, \ell \} = \{ \jmath, m\}$. Then $(A, \rho, t)$ is a balanced oriented quantum algebra if and only if:
\begin{enumerate}
\item[{\rm a)}]
$\rho_{\imath \, \jmath \, \imath \, \jmath} \neq 0$ for all $1 \leq \imath, \jmath \leq n$.
\item[{\rm b)}]
There is a partition $S_1, \ldots, S_r$ of $\{ 1, \ldots, n\}$ into non-empty subsets such that for each $S_\ell$ with $|S_\ell| \geq 2$ the relation $\imath \prec \jmath$ if and only if  $\imath \neq \jmath$ and $\rho_{\imath \, \jmath \, \jmath \, \imath} \neq 0$ determines a well-ordering $\preceq$ on $S_\ell$.
\item[{\rm c)}]
$\rho_{\imath \, \jmath \, \jmath \, \imath} = 0$ unless $\imath, \jmath \in S_\ell$ for some $1 \leq \ell \leq r$.
\item[{\rm d)}]
If  $|S_\ell| \geq 2$ then there exist $b\!c_\ell, x_\ell \in k^\star$ such that
\begin{enumerate}
\item[{\rm i)}]
for all $u \in S_\ell{\setminus}e$ the formula
$$
\omega_u^2 = \left(\frac{\rho_{e \, e \, e \, e}\rho_{u \, u \, u \, u}}{b\!c_\ell}\right)\left(\prod_{e \prec \jmath \prec u} \frac{(\rho_{\jmath \, \jmath \, \jmath \, \jmath})^2}{b\!c_\ell}\right)\omega_e^2
$$
holds, where $e$ is the smallest element in $S_\ell$ with respect to the well-ordering $\preceq$ on $S_\ell$,
\end{enumerate}
and for all $\imath, \jmath \in {\cal S}$ such that $\imath \prec \jmath$:
\begin{enumerate}
\item[{\rm ii)}]
$\rho_{\imath \, \jmath \, \imath \, \jmath}\rho_{\jmath \, \imath \, \jmath \, \imath} = b\!c_\ell$, 
\item[{\rm iii)}]
$\rho_{\imath \, \jmath \, \jmath \, \imath} = x_\ell = \rho_{\imath \, \imath \, \imath \, \imath} - b\!c_{\ell}/ \rho_{\imath \, \imath \, \imath \, \imath}$, and 
\item[{\rm iv)}]
either $\rho_{\imath \, \imath \, \imath \, \imath} = \rho_{\jmath \, \jmath \, \jmath \, \jmath}$ or $\rho_{\imath \, \imath \, \imath \, \imath}\rho_{\jmath \, \jmath \, \jmath \, \jmath} = - b\!c_\ell$. 
\end{enumerate}
\end{enumerate} \qed
\end{TH}
\section{Invariants Associated with $\rho$-Matrices of the Previous Section}\label{SecINvR}
Suppose that $n \geq 2$ and that $(A, \rho, t)$ is the balanced oriented quantum algebra of Theorem \ref{ThmLast}. By Lemma \ref{LEMG} there is an invertible $G \in A = {\rm M}_n(k)$ such that $(A, \rho, t, G)$ is a twist balanced oriented quantum algebra over $k$ and $G$ is unique up to scalar multiple. In this section we determine the resulting link invariants defined in Section \ref{secInv}, where ${\sf tr}$ is the trace function on $A$, in the case of one $\rho$-component.

We continue with the notation of the previous section and will use parts of Theorem \ref{ThmLast} without particular reference. We may assume that $\preceq$ is the usual well-ordering of $\{ 1, \ldots, n\}$ and we continue with the notation $\rho_{\imath \, \imath \, \imath \, \imath} = a_{\imath}$ for all $1 \leq \imath \leq n$. We set $a_1 = a$, $b\!c_1 = b\!c$, $x_1 = x$ and assume that $b\!c$ has a square root in $k$ which we denote by $\sqrt{b\!c}$. Let 
$$
r = \frac{a^2}{b\!c},  \qquad {\bf F} = {\bf Inv}_{A, \, {\sf tr}} \qquad \mbox{and} \qquad q = \frac{a}{\sqrt{b\!c}}.
$$
Without loss of generality we set $G =  \sum_{\imath = 1}^n \omega_\imath^2E_{\imath \, \imath}$. Note that the most general twist balanced oriented quantum algebra structure on $(A, \rho, t)$ is $(A, \rho, t, zG)$, where $z \in k^\star$. If ${\bf Inv}_{A, \, {\sf tr}}, {\bf Inv}'_{A, \, {\sf tr}} : {\SL Link} \la k$ are the functions of Section \ref{KnotLinkINV} defined for $(A, \rho, t, G)$ and $(A, \rho, t, zG)$ respectively, then 
$$
{\bf Inv}'_{A, \, {\sf tr}}({\bf L}) = z^{{\rm Wd}\,{\bf L}_1 + \cdots + {\rm Wd}\,{\bf L}_r}{\bf Inv}_{A, \, {\sf tr}}({\bf L}),
$$
where ${\bf L}_1, \ldots, {\bf L}_r$ are the components of ${\bf L} \in {\SL Link}$.

Let $B = (1/\sqrt{b\!c})(\sum_{\imath, \jmath, \ell, m = 1}^n \rho_{\ell \, \imath \, \jmath \, m} E_{\imath \, \jmath} {\otimes} E_{\ell \, m})$. Then $B$ satisfies the braid equation, is invertible and $B^{-1} = \sqrt{b\!c}(\sum_{\imath, \jmath, \ell, m = 1}^n \rho^{-1}_{\imath \, \ell \, m \, \jmath} E_{\imath \, \jmath} {\otimes} E_{\ell \, m})$. Since $B - B^{-1} = (q - q^{-1})1_A {\otimes} 1_A$, the oriented link invariant ${\bf G}$ defined by 
$$
{\bf G}({\bf L}) = \left(\frac{1}{\sqrt{b\!c}}\right)^{{\rm writhe} \, {\bf L}}{\bf F}({\bf L})
$$
for all ${\bf L} \in {\SL Link}$ satisfies the skein relation
\begin{equation}\label{Eqskein}
{\bf G}_{
\begin{picture}(22,30)(0,0)
\put(0,0){\vector(1,1){22}}
\put(22,0){\line(-1,1){10}}
\put(10,12){\vector(-1,1){10}}
\end{picture} }
\quad - \quad 
{\bf G}_{
\begin{picture}(22,30)(0,0)
\put(12,12){\vector(1,1){10}}
\put(0,0){\line(1,1){10}}
\put(22,0){\vector(-1,1){22}}
\end{picture} }
\quad = \quad 
(q - q^{-1}){\bf G}_{
\begin{picture}(10,30)(0,0)
\put(0,0){\vector(0,1){22}}
\put(10,0){\vector(0,1){22}}
\end{picture} }
\end{equation}
We will show that ${\bf G}$ is determined by the Alexander polynomial in one variable or the regular isotopy HOMFLY polynomial in two variables, up to a writhe and/or rotation factor. Let ${\bf T}_{r, +}$, ${\bf T}_{\ell, -}$, ${\bf T}_{\ell, +}$ and ${\bf T}_{r, -}$ be the oriented $1$--$1$ tangle diagrams depicted by 
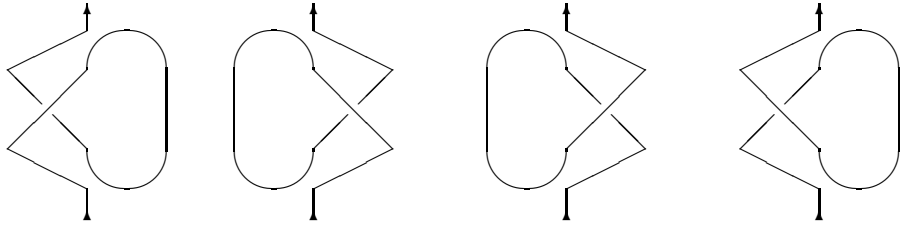
\begin{center}
\begin{picture}(70,100)(-10,-10)
\put(45,15){\oval(30,30)[b]}
\put(45,45){\oval(30,30)[t]}
\put(60,15){\line(0,1){30}}
\put(30,0){\line(0,-1){10}}
\put(30,-8){\vector(0,1){0}}
\put(30,60){\line(0,1){10}}
\put(30,70){\vector(0,1){0}}
%
%
\put(0,15){\line(1,1){30}}
\put(30,15){\line(-1,1){13}}
\put(0,45){\line(1,-1){13}}
%
%
\put(0,15){\line(2,-1){30}}
\put(0,45){\line(2,1){30}}
\end{picture} 
\quad
\begin{picture}(80,100)(-10,-10)
\put(15,15){\oval(30,30)[b]}
\put(15,45){\oval(30,30)[t]}
\put(0,15){\line(0,1){30}}
\put(30,0){\line(0,-1){10}}
\put(30,-8){\vector(0,1){0}}
\put(30,60){\line(0,1){10}}
\put(30,70){\vector(0,1){0}}
%
%
\put(30,15){\line(1,1){13}}
\put(60,15){\line(-1,1){30}}
\put(60,45){\line(-1,-1){13}}
%
%
\put(30,0){\line(2,1){30}}
\put(30,60){\line(2,-1){30}}
\end{picture} 
\quad
\begin{picture}(80,100)(-10,-10)
\put(15,15){\oval(30,30)[b]}
\put(15,45){\oval(30,30)[t]}
\put(0,15){\line(0,1){30}}
\put(30,0){\line(0,-1){10}}
\put(30,-8){\vector(0,1){0}}
\put(30,60){\line(0,1){10}}
\put(30,70){\vector(0,1){0}}
%
%
\put(30,15){\line(1,1){30}}
\put(60,15){\line(-1,1){13}}
\put(30,45){\line(1,-1){13}}
%
%
\put(30,0){\line(2,1){30}}
\put(30,60){\line(2,-1){30}}
\end{picture} 
\quad
\begin{picture}(70,100)(-10,-10)
\put(45,15){\oval(30,30)[b]}
\put(45,45){\oval(30,30)[t]}
\put(60,15){\line(0,1){30}}
\put(30,0){\line(0,-1){10}}
\put(30,-8){\vector(0,1){0}}
\put(30,60){\line(0,1){10}}
\put(30,70){\vector(0,1){0}}
%
%
\put(0,15){\line(1,1){13}}
\put(30,45){\line(-1,-1){13}}
\put(30,15){\line(-1,1){30}}
%
%
\put(0,15){\line(2,-1){30}}
\put(0,45){\line(2,1){30}}
\end{picture} 
\end{center}
respectively. For $m \geq 0$ let ${\sf C}_{r, +}(m)$, ${\sf C}_{r, -}(m)$ be the oriented diagrams 
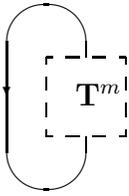
\begin{picture}(60,80)(-5,-25)
\put(0,-5){\line(0,1){40}}
\put(15,-5){\oval(30, 30)[b]}
\put(15,35){\oval(30, 30)[t]}
\put(30,0){\line(0,-1){5}}
\put(30,30){\line(0,1){5}}
\put(15,0){\dashbox{5}(30,30)}
\put(26,12){${\bf T}^m$}
\put(0,15){\vector(0, -1){0}}
\end{picture}
where ${\bf T} = {\bf T}_{r, +}$, ${\bf T}_{r, -}$ respectively, and let ${\sf C}_{\ell, +}(m)$, ${\sf C}_{\ell, -}(m)$ be the oriented diagram
\begin{picture}(60,80)(-5,-25)
\put(45,-5){\line(0,1){40}}
\put(30,-5){\oval(30, 30)[b]}
\put(30,35){\oval(30, 30)[t]}
\put(15,0){\line(0,-1){5}}
\put(15,30){\line(0,1){5}}
\put(0,0){\dashbox{5}(30,30)}
\put(11,12){${\bf T}^m$}
\put(45,15){\vector(0, -1){0}}
\end{picture}
where ${\bf T} = {\bf T}_{\ell, +}$, ${\bf T}_{\ell, -}$ respectively. By convention ${\sf C}_{\ell, \pm}(0)$ is a circle with clockwise orientation and ${\sf C}_{r, \pm}(0)$ is a circle with counter clockwise orientation.  To determine ${\bf F}$ we compute ${\bf F}$ on these four types of composite diagrams and compute ${\rm Tr}(G)$, ${\rm Tr}(G^{-1})$.

For $0 \leq \ell \leq m \leq n$ let 
$$
\eta_{+} (\ell {:} m) = |\{ \imath \, | \, \ell < \imath \leq m, a_\imath = a\}|
$$
and 
$$
\eta_{-} (\ell {:} m) = |\{ \imath \, | \, \ell < \imath \leq m, a_\imath \neq a\}|.
$$
Now for each $1 \leq \imath \leq n$ either $a_\imath = a$ or $a_\imath = -b\!c/a$. Thus $a_\imath^2/b\!c = r$ if $a_\imath = a$, and $a_\imath^2/b\!c = r^{-1}$ if $a_\imath \neq a$. Consequently 
\begin{equation}\label{EqwiSquared}
\omega_\imath^2 = - (-r)^{-\delta^{a_\imath}_a}r^{\eta_{+} (0 {:} \imath) - \eta_{-} (0 {:} \imath)}
\end{equation}
for all $1 \leq \imath \leq n$. Set $\omega_\imath (r) = \omega_\imath$, and otherwise for $x \in k^\star$ choose $\omega_\imath (x) \in k$ so that $\omega_\imath (x)^2 = - (-x)^{-\delta^{a_\imath}_a}x^{\eta_{+} (0 {:} \imath) - \eta_{-} (0 {:} \imath)}$. For all $0 \leq \ell \leq n$ and $x \neq 0, 1$ we have the formula
\begin{equation}\label{EqSum1x}
\sum_{\jmath = 1}^\ell \omega_\jmath (x)^2 = \frac{1 - x^{\eta_{+} (0 {:} \ell) - \eta_{-} (0 {:} \ell)}}{1 - x},
\end{equation}
and thus 
\begin{equation}\label{EqSum2x}
\sum_{\jmath = \ell + 1}^n \omega_\jmath (x)^2 = x^{\eta_{+} (0 {:} \ell) - \eta_{-} (0 {:} \ell)}\left( \frac{1 - x^{\eta_{+} (\ell {:} n) - \eta_{-} (\ell {:} n)}}{1 - x}\right).
\end{equation}
To establish (\ref{EqSum1x}) observe that $\omega_{\imath} (x)^2 + \omega_{\imath + 1} (x)^2 = 0$ whenever $a_\imath$ and $a_{\imath + 1}$ are different. Thus for the purpose of establishing (\ref{EqSum1x}) we may assume that the $a_\imath$'s in the list $a_1 , \ldots, a_\ell$ which are equal to $a$ precede those $a_\imath$'s which are not equal to $a$. 

Observe that $\omega_\imath (1)^2 = (-1)^{\delta^{a_\imath}_a + 1}$ and $\omega_\imath (xy)^2 = (-1)^{\delta^{a_\imath}_a + 1} \omega_\imath (x)^2 \omega_\imath (y)^2$ for all $x, y \in k^\star$. The transparent relation
\begin{equation}\label{Eqw2andw2Inv}
\omega_\imath (x^{-1})^2 = \omega_\imath (x)^{-2}
\end{equation}
is also a consequence of the preceding properties which the $\omega_\imath (x)$'s satisfy. Let 
$$
\hbar = \frac{r^{\eta_{+} (0 {:} n) - \eta_{-} (0 {:} n)}}{r}.
$$
By (\ref{EqSum1x}) and (\ref{Eqw2andw2Inv}) we see that 
$$
{\rm Tr}(G) = \frac{1 - r^{\eta_{+} (0 {:} n) - \eta_{-} (0 {:} n)}}{1 - r}, \quad 
{\rm Tr}(G^{-1}) = \frac{1 - r^{-(\eta_{+} (0 {:} n) - \eta_{-} (0 {:} n))}}{1 - r^{-1}}
$$
and consequently 
$$
\frac{{\rm Tr}(G)}{{\rm Tr}(G^{-1})} = \hbar
$$
when ${\rm Tr}(G^{-1}) \neq 0$. Note that ${\rm Tr}(G^{-1}) \neq 0$ if and only if  ${\rm Tr}(G) \neq 0$. Using (\ref{EqSum2x}) and (\ref{EqSum1x}) one can show that 
$$
{\bf w}_A({\bf T}_{r, +}) = a(\sum_{\imath = 1}^n (-r)^{\delta^{a_\imath}_a - 1}r^{\eta_{+} (\imath {:} n) - \eta_{-} (\imath {:} n)}E_{\imath \, \imath} ), \;\; {\bf F}({\sf C}_{r, +}(m)) = (a\hbar )^m {\rm Tr}(G^{-1}),
$$
$$
{\bf w}_A({\bf T}_{\ell, -}) = \frac{1}{a}(\sum_{\imath = 1}^n (-r^{-1})^{\delta^{a_\imath}_a - 1}r^{-(\eta_{+} (\imath {:} n) - \eta_{-} (\imath {:} n))}E_{\imath \, \imath}), \;\; {\bf F}({\sf C}_{\ell, -}(m)) = (a\hbar )^{-m}{\rm Tr}(G),
$$
$$
{\bf w}_A({\bf T}_{\ell, +}) = a(\sum_{\imath = 1}^n r^{-1}(-r^{-1})^{\delta^{a_\imath}_a - 1}r^{\eta_{+} (0 {:} \imath ) - \eta_{-} (0 {:} \imath )}E_{\imath \, \imath}), \;\; {\bf F}({\sf C}_{\ell, +}(m)) = a^m{\rm Tr}(G),
$$
and
$$
{\bf w}_A({\bf T}_{r, -}) = \frac{1}{a}(\sum_{\imath = 1}^n r(-r)^{\delta^{a_\imath}_a - 1}r^{-(\eta_{+} (0 {:} \imath ) - \eta_{-} (0 {:} \imath ))}E_{\imath \, \imath}), \;\; {\bf F}({\sf C}_{r, -}(m)) = a^{-m}{\rm Tr}(G^{-1})
$$
for all $m \geq 0$. 

Suppose that ${\rm Tr}(G) = 0$. Then ${\bf G}({\bf L}) = \nabla_{\bf L}(q - q^{-1})$ for all ${\bf L} \in {\SL Link}$, where $\nabla_{\bf L}(z)$ is the one variable Alexander polynomial. See 
\cite[page 174]{KNOTS}
for example. Thus 
$$
{\bf F}({\bf L}) = a^{{\rm writhe} \, {\bf L}}\left( q^{- {\rm writhe} \, {\bf L}} \nabla_{\bf L}(q - q^{-1})\right)
$$
for all ${\bf L} \in {\SL Link}$ when ${\rm Tr}(G) = 0$.

Now suppose that ${\rm Tr}(G) \neq 0$. Then ${\rm Tr}(G^{-1}) \neq 0$. Set 
$$
\rho = \frac{q}{q^{\eta_{+} (0 {:} n) - \eta_{-} (0 {:} n)}} \quad \mbox{and} \quad \kappa = \rho {\rm Tr}(G) = \rho^{-1} {\rm Tr}(G^{-1}).
$$
At this point it is not hard to see that ${\bf F}({\sf C}) = \kappa (a\rho^{-1})^{{\rm writhe} \, {\sf C}}\rho^{- {\rm Wd} \, {\sf C}}$ for all ${\sf C} = {\sf C}_{\ell, \pm}(m),  {\sf C}_{r, \pm}(m)$, where $m \geq 0$. Using (\ref{Eqskein}) we see that 
$$
H_{\bf L}(q^{\eta_{+} (0 {:} n) - \eta_{-} (0 {:} n)}, q - q^{-1}) = \frac{(1/\sqrt{b\!c})^{{\rm writhe} \, {\bf L}}\rho^{{\rm Wd} \, {\bf L}}}{\kappa} {\bf F}({\bf L})
$$
for all ${\bf L} \in {\SL Link}$, where $H_{\bf L}(\alpha, z)$ is the two variable regular isotopy HOMFLY polynomial. See 
\cite[page 54]{KNOTS}
for example. Thus 
$$
{\bf F}({\bf L}) = a^{{\rm writhe} \, {\bf L}}\left( \kappa q^{- {\rm writhe} \, {\bf L}}\rho^{- {\rm Wd} \, {\bf L}} H_{\bf L}( q^{\eta_{+} (0 {:} n) - \eta_{-} (0 {:} n)}, q - q^{-1})\right)
$$
for all ${\bf L} \in {\SL Link}$ when ${\rm Tr} (G) \neq 0$.

Suppose that $a$ and $\sqrt{b\!c}$ are independent indeterminates over the prime field of $k$. Note that ${\bf F}({\bf L})$ is a polynomial in $a$ and $b\!c$ and is a homogeneous Laurent polynomial in $a$ and $\sqrt{b\!c}$ of degree ${\rm writhe} \, {\bf L}$. By contrast the formulas for ${\bf F}$ deduced in this section are in terms of $a$ and $\sqrt{b\!c}$.

\end{document}